\documentclass[journal,12pt,onecolumn,]{IEEEtran}
\usepackage{cite}
\usepackage{color}
\usepackage[pdftex]{graphicx}
\graphicspath{{fig/}{jpeg/}}
\usepackage[cmex10]{amsmath}
\usepackage{amssymb,bbm}
\usepackage{algorithm,algorithmic,multirow,hhline}
\usepackage{amsmath,amssymb,lipsum}
\usepackage{hyperref}
\usepackage{comment}
\usepackage{graphicx}
\usepackage[font=small]{caption}
\usepackage{subcaption}
\usepackage{mathtools}

\input{mysymbol.sty}
\usepackage{needspace}

\usepackage{theorem}

\newtheorem{thm}{Theorem}
\newtheorem{lemma}{Lemma}
\newtheorem{proposition}{Proposition}
\newtheorem{corollary}{Corollary}

\theoremstyle{definition}
\newtheorem{assumption}{Assumption}
\newtheorem{remark}{Remark}

\def \QED {$\blacksquare$}
\newcommand{\closure}[2][3]{{}\mkern#1mu\overline{\mkern-#1mu#2}}

 \newcommand{\INDSTATE}[1][1]{\STATE\hspace{3mm}}
\newcommand{\INDSTATED}[1][1]{\STATE\hspace{6mm}}

\title{Decentralized Online Learning with Kernels}
\author{Alec Koppel$^\S$, Santiago Paternain$^{\star}$, C\'{e}dric Richard$^\dagger$ and Alejandro Ribeiro$^\star$
\thanks{This work in this paper is supported by NSF CCF-1017454, NSF CCF-0952867, ONR N00014-12-1-0997, ARL MAST CTA, and ASEE SMART.Part of the results in this paper appeared in \cite{koppel2017multipolk}.}
\thanks{ \noindent$^\S$Computational and Information Sciences Directorate, U.S. Army Research Laboratory, Adelphi, MD, 20783. Email: alec.e.koppel.civ@mail.mil}
\thanks{ \noindent$^\star$Department of ESE, University of Pennsylvania, 200 South 33rd Street, Philadelphia, PA 19104. Email: \{spater, aribeiro\}@seas.upenn.edu}
\thanks{\noindent $^{\dagger}$Laboratory Lagrange - UMR CNRS 7293, Observatory of the French Riviera	 
University of Nice Sophia-Antipolis, Nice, France, 06108} 
}

\begin{document}
\maketitle

\begin{abstract}
We consider multi-agent stochastic optimization problems over reproducing kernel Hilbert spaces (RKHS). In this setting, a network of interconnected agents aims to learn decision functions, i.e., nonlinear statistical models, that are optimal in terms of a global convex functional that aggregates data across the network, with only access to locally and sequentially observed samples. We propose solving this problem by allowing each agent to learn a local regression function while enforcing consensus constraints. 
We use a penalized variant of functional stochastic gradient descent operating simultaneously with low-dimensional subspace projections. These subspaces are constructed greedily by applying orthogonal matching pursuit to the sequence of kernel dictionaries and weights. By tuning the projection-induced bias, we propose an algorithm that allows for each individual agent to learn, based upon its locally observed data stream and message passing with its neighbors only, a regression function that is close to the globally optimal regression function.
That is, we establish that with constant step-size selections agents' functions converge to a neighborhood of the globally optimal one while satisfying the consensus constraints as the penalty parameter is increased. Moreover,  the complexity of the learned regression functions is guaranteed to remain finite.
On both multi-class kernel logistic regression and multi-class kernel support vector classification with data generated from class-dependent Gaussian mixture models, we observe stable function estimation and state of the art performance for distributed online multi-class classification. Experiments on the Brodatz textures further substantiate the empirical validity of this approach. 

\end{abstract}

\section{Introduction}\label{sec:intro}

We consider decentralized online optimization problems: a network $\mathcal{G}=(V,\mathcal{E})$ of agents aims to minimize a global objective that is a sum of local convex objectives available only to each node. The problem is online and distributed because data samples upon which the local objectives depend are sequentially and locally observed by each agent. In this setting, agents aim to make inferences as well as one which has access to all data at a centralized location in advance. Instead of assuming agents seek a common parameter vector $\bbw\in\reals^p$, we focus on the case where agents seek to learn a common \emph{decision function} $f(\bbx)$ that belong to a reproducing kernel Hilbert space (RKHS). Such functions represent, e.g., nonlinear statistical models \cite{anthony2009neural} or trajectories in a continuous space \cite{marinhofunctional}. 
Learning in multi-agent settings arises predominately in two technological settings: industrial-scale machine learning, where optimizing statistical model parameters is decentralized across a parallel processing architecture to attain computational speedup; and networked intelligent systems such as sensor networks \cite{kozick2004source}, multi-robot teams \cite{cKoppelEtal16a,schwager2017multi}, and Internet of Things \cite{liu2012algorithm,ghosh2015pricing}. In the later setting, decentralized processing justified as opposed to using a fusion center when the communication cost of centralization exceeds the cost of distributed information protocols. This is true of multi-agent systems with streaming data considered here.

Efforts to develop optimization tools for multi-agent online learning have thus far been restricted to the case where each agent learns a linear statistical model \cite{Koppel2015a} or a task-driven dictionary \cite{KoppelEtal16a} that is as good as one with data aggregated across the network. However, these efforts exclude the state of the art tools for statistical learning based on nonlinear interpolators: namely, kernel methods \cite{slavakis2013online,li2014kernel} and neural networks \cite{haykin1994neural,krizhevsky2012imagenet}. We note that instabilities associated with non-convexity which are only a minor issue in centralized settings \cite{Mairal2012} become both theoretically and empirically difficult to overcome in settings with consensus constraints \cite{KoppelEtal16a}, and therefore efforts to extend neural network learning to multi-agent online learning likely suffer the same drawbacks.\footnote{In general, globally convergent decentralized online training of neural networks is an open problem, whose solution requires fundamentally new approaches to stochastic global optimization.} Therefore, we focus on extending kernel methods to decentralized online settings, motivated both by its advantageous empirical performance, as well as the theoretical and practical benefits of the fact that the optimization problem defined by their training is convex. This stochastic convex problem, however, is defined over an infinite dimensional space, and therefore it is not enough to solve the optimization problem, but one must also solve it in an optimally sparse way. Doing so in multi-agent settings is the goal of this work.

To contextualize our solution methodology, consider centralized vector-valued stochastic convex programming, which has classically been solved with stochastic gradient descent (SGD) \cite{Robbins1951}. SGD involves descending along the negative of the stochastic gradient rather than the true gradient to avoid the fact that computing the gradient of the average objective has complexity comparable to the training sample size, which could be infinite. 
In contrast, the setting considered in this work is a stochastic program defined over a function space, which is in general an intractable variational inference problem. However, when the function space is a RKHS \cite{kimeldorf1971some}, the Representer Theorem allows us to transform a search over an infinite space into one over a set of weights and data samples \cite{scholkopfgeneralized}. Unfortunately, the feasible set of the resulting problem has complexity comparable to the sample size $N$, and thus is intractable for $N\rightarrow \infty$ \cite{norkin2009stochastic}. Compounding this problem is that the storage required to construct the functional generalization of SGD is comparable to the iteration index of the algorithm, which is untenable for online settings.

Efforts to mitigate the complexity of the function representation (``the curse of kernelization") have been previously developed. These combine functional extensions of stochastic gradient method with compressions of the function parameterization independently of the optimization problem to which they are applied \cite{1315946, liu2008kernel, Kivinen2004, export:78226, Zhu2005} or approximate the kernel during training \cite{dai2014scalable,le2016nonparametric,le2016dual,lu2016large,calandriello2017second}, and at best converge on average.
In contrast, a method was recently proposed that combines greedily constructed \cite{Pati1993} sparse subspace projections with functional stochastic gradient method and guarantees exact convergence to the minimizer of the  average risk functional. This technique, called parsimonious online learning with kernels (POLK), tailors the parameterization compression to preserve the descent properties of the underlying RKHS-valued stochastic process  \cite{POLK}, and inspires the approach considered here.

In this work, we extend the ideas in \cite{POLK} to multi-agent settings. Multiple tools from distributed optimization may be used to do so; however, we note that the Representer Theorem \cite{scholkopfgeneralized} has not been established for general stochastic saddle point problems in RKHSs. Therefore, we adopt an approximate primal-only approach based on penalty methods \cite{Johansson2008,Ram2010}, which in decentralized optimization is known as distributed gradient descent (DGD). Using functional stochastic extensions of DGD, together with the greedy Hilbert subspace projections designed in POLK, we develop a method such that each agent, through its local data stream and message passing with only its neighbors, learns a memory-efficient approximation to the globally optimal regression function with probability $1$. Such global stability guarantees are in contrast to specialized results for multi-agent kernel learning \cite{nguyen2005nonparametric,forero2010consensus} and alternative distributed online nonlinear function estimation methods such as dictionary learning \cite{Mairal2012,Chainais2013,KoppelEtal16a} or neural networks \cite{krizhevsky2012imagenet}, which suffer from instability due to the non-convexity of the optimization problem their training defines.

The result of the paper is organized as follows. In Section \ref{sec:prob} we clarify the problem setting of stochastic programming in RKHSs in the centralized and decentralized case. In Section \ref{sec:algorithm}, we propose a new penalty functional that permits deriving a decentralized online method for kernel regression without any complexity bottleneck by making use of functional stochastic gradient method (Section \ref{subsec:fsg}) combined with greedy subspace projections (Section \ref{subsec:proj}). In Section \ref{sec:convergence} we present our main theoretical results, which establishes that the function sequence of each agent generated by the proposed technique converges to a neighborhood of the globally optimal function with probability $1$. In Section \ref{sec:simulations}, we present numerical examples of decentralized online multi-class kernel logistic regression and kernel support vector machines with data generated from Gaussian mixtures, and observe a state of the art trade-off between Lyapunov stability and statistical accuracy. We then apply the resulting method to the benchmark Brodatz texture dataset \cite{Brodatz1966} and observe state of the art decentralized online multi-class classification performance.

\section{Problem Formulation}\label{sec:prob}
\subsection{Decentralized Functional Stochastic Programming}\label{subsec:opt}

Consider the problem of expected risk minimization, where the goal is to learn a regressor that minimizes a loss function quantifying the merit of a statistical model averaged over a data set. We focus on the case when the number of training examples $N$ is very large or infinite. In this work, input-output examples, $(\bbx_n, \bby_n)$, are i.i.d. realizations drawn from a stationary joint distribution over the random pair $(\bbx, \bby) \in \ccalX \times \ccalY$, where $\ccalX \subset \reals^p$ and $\ccalY \subset \reals$. Here, we consider finding regressors that are not vector valued parameters, but rather functions $\tilde{f} \in \ccalH$ in a hypothesized function class $\ccalH$, which allows for learning nonlinear statistical models rather than generalized linear models that rarely achieve satisfactory statistical error rates in practice \cite{mukherjee1996automatic,li2014kernel}. The merit of the function $\tilde{f}$ is evaluated by the convex loss function $\ell:\ccalH \times \ccalX \times \ccalY \rightarrow \reals$ that quantifies the merit of the estimator $\tilde{f}(\tbx)$ evaluated at feature vector $\tbx$. This loss is averaged over all possible training examples to define the statistical loss $\tilde{L}(\tilde{f}) : = \mbE_{\bbx, \bby}{[ \ell(\tilde{f}(\bbx), y)}]$, which we combine with  a Tikhonov regularizer to construct the regularized loss $\tilde{R}(\tilde{f}) : = \argmin_{\tilde{f}\in \ccalH} \tilde{L}(\tilde{f}) +(\lambda/2)\|\tilde{f}\|^2_{\ccalH}$ \cite{shalev2010learnability, evgeniou2000regularization}.
We then define the optimal function as\vspace{-2mm}
\begin{align}\label{eq:kernel_stoch_opt}
\!\!\!\tilde{f}^*\!=\!\argmin_{\tilde{f} \in \ccalH} \tilde{R}(\tilde{f}) : =\! \argmin_{\tilde{f} \in \ccalH}\mbE_{\tbx, \tilde{y}}\!\Big[ \ell(\tilde{f}\big(\tbx), \tilde{y}\big)\Big] \!+\!\frac{\lambda}{2}\|\tilde{f} \|^2_{\ccalH} 
\end{align}
In this work, we focus on extensions of the formulation in \eqref{eq:kernel_stoch_opt} to the case where data is scattered across an interconnected network that represents, for instance, robotic teams \cite{KoppelEtal16a}, communication systems \cite{ribeiro2010ergodic}, or sensor networks \cite{kozick2004source}. To do so, we define a symmetric, connected, and directed network $\ccalG = (\ccalV, \ccalE)$ with $|\ccalV|=V$ nodes and $|\ccalE|=E$ edges and denote as $n_i:=\{j:(i,j)\in\ccalE\}$ the neighborhood of agent $i$. For simplicity we assume that the number of edges $E$ is even. Each agent $i\in \ccalV$ observes a local data sequence as realizations $(\bbx_{i,n}, y_{i,n})$ from random pair $(\bbx_i, y_i) \in \ccalX \times \ccalY$ and seeks to learn a common globally optimal regression function $f$. This setting may be mathematically captured by associating to each node $i$ a convex loss functional $\ell_i:\ccalH \times \ccalX \times \ccalY \rightarrow \reals$ that quantifies the merit of the estimator $f_i(\bbx_i)$ evaluated at feature vector $\bbx_i$, and defining the goal for each node as the minimization of the common global loss
\begin{align}\label{eq:kernel_stoch_opt_global}
f^*= \argmin_{f \in \ccalH}\sum_{i\in\ccalV}\left(\mbE_{\bbx_i,y_i}\Big[ \ell_i(f\big(\bbx_i), y_i\big)\Big] +\frac{\lambda}{2}\|f \|^2_{\ccalH} \right)
\end{align}
Observe that this global loss is a network-wide average (scaled by $V$) of all local losses, and therefore the minimizers of \eqref{eq:kernel_stoch_opt} and \eqref{eq:kernel_stoch_opt_global} coincide when $(\bbx_i, y_i)$ have a common joint distribution for each $i$. However, in multi-agent optimization, this is not generally the case, thus when selecting a regression function $f$ with only local data, different agents will learn a different decision function $f_i^*$ that it is not optimal as compared to one selected in a centralized manner, i.e., with the data gathered by all agents. To overcome this limitation we allow message passing between agents and we impose a consensus constraint on the regression function among neighbors  $f_i=f_j \;, \ (i,j)\in \ccalE$. Thus we consider the nonparametric decentralized stochastic program:
\begin{align}\label{eq:main_prob}
f^*=& \ \ \argmin_{ \{f_i\} \subset \ccalH} \qquad \sum_{i\in\ccalV}\left(\mbE_{\bbx_i, y_i}\Big[ \ell_i(f_i\big(\bbx), y_i\big)\Big] +\frac{\lambda}{2}\|f_i \|^2_{\ccalH} \right) \nonumber \\
&\quad \text{such that}\quad \ f_i=f_j\;,  (i,j)\!\in \ccalE
\end{align}
For further define the product Hilbert space $\ccalH^V$ of functions aggregated over the network whose elements are stacked functions $f(\cdot) = [f_1(\cdot) ; \cdots ; f_V(\cdot)]$ that yield vectors of length $V$ when evaluated at local random vectors $f(\bbx) = [f_1(\bbx_1) ; \cdots ; f_V(\bbx_V)] \in \reals^V$. Moreover, define the stacked random vectors $\bbx = [\bbx_1 ; \cdots ; \bbx_V] \in \ccalX^V \subset \reals^{Vp}$ and $\bby = [ y_1 ; \cdots y_V] \in \reals^V$ that represents $V$ labels or physical measurements, for instance. 

 The goal of this paper is to develop an algorithm to solve \eqref{eq:main_prob} in distributed online settings where nodes do not know the distribution of the random pair $(\bbx_i, y_i)$ but observe local independent training examples $(\bbx_{i,n}, y_{i,n})$ sequentially. 
\subsection{Function Estimation in Reproducing Kernel Hilbert Spaces}\label{subsec:kernels}
The optimization problem in \eqref{eq:kernel_stoch_opt}, and hence \eqref{eq:main_prob}, is intractable in general, since it defines a variational inference problem integrated over the unknown joint distribution $\mathbb{P}(\bbx,y)$.
 However, when $\ccalH$ is equipped with a \emph{reproducing kernel} $\kappa : \ccalX \times \ccalX \rightarrow \reals$ (see \cite{MullerAFPT13,li2014kernel}), a function estimation problem of the form \eqref{eq:kernel_stoch_opt} may be reduced to a parametric form via the Representer Theorem \cite{wheeden1977measure,norkin2009stochastic}.
Thus, we restrict the Hilbert space in Section \ref{subsec:opt} to be one equipped with a kernel $\kappa$ that satisfies for all functions $\tilde{f}:\ccalX\rightarrow \reals$ in $\ccalH$:
\begin{align} \label{eq:rkhs_properties}
& (i)  \  \langle \tilde{f} , \kappa(\bbx_i, \cdot)) \rangle _{\ccalH} = \tilde{f}(\bbx_i) ,\quad (ii)  \ \ccalH = \closure{\text{span}\{ \kappa(\bbx_i , \cdot) \}}  
\end{align}
for all $ \bbx_i \in \ccalX$. Here $\langle \cdot , \cdot \rangle_{\ccalH}$ denotes the Hilbert inner product for $\ccalH$. Further assume that the kernel is positive semidefinite, i.e. $\kappa(\bbx_i, \bbx_i') \geq 0$ for all $\bbx_i, \bbx_i' \in \ccalX$. Function spaces of this type are called reproducing kernel Hilbert spaces (RKHS).

In \eqref{eq:rkhs_properties}, property (i) is the reproducing property (via Riesz Representation Theorem \cite{wheeden1977measure}).  Replacing $\tilde{f}$ by $\kappa(\bbx_i' , \cdot) $  in \eqref{eq:rkhs_properties} (i) yields $ \langle \kappa(\bbx_i', \cdot) , \kappa(\bbx_i, \cdot) \rangle_{\ccalH} = \kappa(\bbx_i, \bbx_i')$ which is the origin of the term ``reproducing kernel."  This property induces a nonlinear transformation of the input space $\ccalX$: denote by $\phi(\cdot)$ a nonlinear map of the feature space that assigns to each $\bbx_i$ the kernel function $\kappa(\cdot, \bbx_i)$. The reproducing property yields that the inner product of the image of distinct feature vectors $\bbx_i$ and $\bbx_i'$ under the map $\phi$ requires only kernel evaluations: $\langle \phi(\bbx_i), \phi(\bbx_i') \rangle_{\ccalH} =\kappa(\bbx_i, \bbx_i')$ (the 'kernel trick').

Moreover, property \eqref{eq:rkhs_properties} (ii) states that functions $\tilde{f}\in \ccalH$ may be written as a linear combination of kernel evaluations. For kernelized and regularized empirical risk minimization (ERM), the Representer Theorem \cite{kimeldorf1971some,scholkopfgeneralized} establishes that the optimal $\tilde{f}$ in  hypothesized function class $\ccalH$ admit an expansion in terms of kernel evaluations \emph{only} over training examples\vspace{-2mm}
%
%
\begin{equation}\label{eq:kernel_expansion}
\tilde{f}(\bbx_i) = \sum_{n=1}^N w_{i,n} \kappa(\bbx_{i,n}, \bbx_i)\; ,
\end{equation}
where $\bbw_i = [w_{i,1}, \cdots, w_{i,N}]^T \in \reals^N$ denotes a set of weights. The upper index $N$ in \eqref{eq:kernel_expansion} is referred to as the model order, and for ERM the model order and training sample size are equal. Common choices $\kappa$ include the polynomial and radial basis kernels, i.e., $\kappa(\bbx_i,\bbx_i') = \left(\bbx_i^T\bbx_i'+b\right)^d $ and $\kappa(\bbx_i,\bbx_i') = \exp\{ -{\lVert \bbx_i - \bbx_i' \rVert_2^2}/{2d^2} \}$, respectively, where $\bbx_i, \bbx_i' \in \ccalX$.

Suppose, for the moment, that we have access to $N$ i.i.d. realizations of the random pairs $(\bbx_i, y_i)$ for each agent $i$ such that the expectation in \eqref{eq:main_prob} is computable, and we further ignore the consensus constraint. Then the objective in \eqref{eq:main_prob} becomes:
\begin{align}\label{eq:kernel_batch_opt}
f^*&=\argmin_{f \in \ccalH^V} \frac{1}{N}\sum_{n=1}^N \sum_{i\in \ccalV} \ell(f_i(\bbx_{i,n}), y_{i,n}) +\frac{\lambda}{2}\|f_i \|^2_{\ccalH}\; 
\end{align}
Then, by substituting the Representer Theorem [cf. \eqref{eq:kernel_expansion}] into \eqref{eq:main_prob}, we obtain that optimizing in $\ccalH^V$ reduces to optimizing over the set of $NV$ weights:
\begin{align}\label{eq:kernel_batch_opt2}
\!\!\!f^*\!\!=\!\!\argmin_{\{\!\bbw_i \!\} \in \reals^N}\!\!\frac{1}{N}\!\!\!\sum_{n=1}^N \! \sum_{i\in\ccalV} \!\! \ell_i(\!\bbw_i^T\!\! \boldsymbol{\kappa}_{\bbX_i}\!(\bbx_{i,n}\!)\!,\! y_{i,n}\!) \!+\!\frac{\lambda}{2}\! \bbw_i^T \bbK_{\!\bbX_i,\bbX_i\!} \!\bbw_i  ,
%
\end{align}
where we have defined the Gram (or kernel) matrix $\bbK_{\bbX_i,\bbX_i}\in \reals^{N\times N}$, with entries given by the kernel evaluations between $\bbx_{i,m}$ and $\bbx_{i,n}$ as $[\bbK_{\bbX_i, \bbX_i}]_{m,n} =\kappa(\bbx_{i,m}, \bbx_{i,n})$.  We further define the vector of kernel evaluations $\boldsymbol{\kappa}_{\bbX_i}(\cdot) = [\kappa(\bbx_{i,1},\cdot) \ldots \kappa(\bbx_{i,N},\cdot)]^T$, which are related to the kernel matrix as $\bbK_{\bbX_i,\bbX_i} = [\boldsymbol{\kappa}_{\bbX_i}(\bbx_{i,1}) \ldots \boldsymbol{\kappa}_{\bbX_i}(\bbx_{i,N})]$. The dictionary of training points associated with the kernel matrix is defined as  $\bbX_i = [\bbx_{i,1},\ \ldots\ ,\bbx_{i,N}]$.

By exploiting the Representer Theorem, we transform a nonparametric infinite dimensional optimization problem in $\ccalH^V$ \eqref{eq:kernel_batch_opt} into a finite $NV$-dimensional parametric problem \eqref{eq:kernel_batch_opt2}. Thus, for empirical risk minimization, the RKHS provides a principled framework to solve nonparametric regression problems as a search over $\reals^{VN}$ for an optimal set of coefficients. 

However, to solve problems of the form \eqref{eq:kernel_batch_opt} when training examples $(\bbx_{i,n}, y_{i,n})$ become sequentially available or their total number $N$ is not finite, the objective in \eqref{eq:kernel_batch_opt} becomes an expectation over random pairs $(\bbx_{i}, y_{i})$ as \cite{slavakis2013online}
\begin{align}\label{eq:kernel_stoch_opt2}
f^* &= \! \!\!\!\argmin_{\bbw_i\in \reals^{\ccalI}, \{\bbx_{i,n}\}_{n\in\ccalI}}\sum_{i \in \ccalV}\mbE_{\bbx_i, y_i}{[ \ell_i(\sum_{n\in\ccalI} w_{i,n} \kappa(\bbx_{i,n}, \bbx_i) , y_i)}]  \nonumber \\
&\qquad \qquad \qquad+\frac{\lambda}{2}\|\!\! \!\sum_{n,m\in \ccalI } \!\!w_{i,n} w_{i,m} \kappa(\bbx_{i,m}, \bbx_{i,n})  \|^2_{\ccalH}  \; , 
\end{align}
where we substitute the Representer Theorem generalized to the infinite sample-size case established in \cite{norkin2009stochastic} into the objective \eqref{eq:main_prob}  with $\ccalI$ as some countably infinite indexing set. That is, as the data sample size $N\rightarrow \infty$, the representation of $f_i$ becomes infinite as well. Thus, our goal is to solve \eqref{eq:kernel_stoch_opt2} in an approximate manner such that each $f_i$ admits a finite representation near $f_i^*$, while satisfying the consensus constraints $f_i = f_j $ for $(i,j)\in \ccalE$ (which were omitted for the sake of discussion between \eqref{eq:kernel_batch_opt} - \eqref{eq:kernel_stoch_opt2}).

\section{Algorithm Development}\label{sec:algorithm}

We turn to developing an online iterative and decentralized solution to solving \eqref{eq:main_prob} when the functions $\{f_i \}_{i \in \ccalV}$ are elements of a RKHS, as detailed in Section \ref{subsec:kernels}. To exploit the properties of this function space, we require the applicability of the Representer Theorem [cf. \eqref{eq:kernel_expansion}], but this result holds for any regularized minimization problem with a convex functional. Thus, we may address the consensus constraint $f_i=f_j\; , (i,j)\in\ccalE$ in \eqref{eq:main_prob} by enforcing approximate consensus on estimates $f_i(\bbx_i) = f_j(\bbx_i)$ in expectation. This specification may be met by introducing the penalty functional
\begin{align}\label{eq:penalty_function}
\psi_c(f)\! &=  \!\!\sum_{i\in\ccalV}\!\!\Big(\!\mbE_{\bbx_i, \bby_i}\! \Big[ \ell_i(f_i\big(\bbx_i), y_i\big)\!\Big]\! \!+\!\frac{\lambda}{2}\|f_i \|^2_{\ccalH} \! \nonumber \\
&\qquad\qquad\! \! + \! \frac{c}{2}  \!\!\sum_{j\in n_i}\!  \mathbb{E}_{\bbx_i}\left\{ [f_i(\bbx_i) \! - \! f_j(\bbx_i)]^2 \right\}\Big)
\end{align}
The reasoning for the definition \eqref{eq:penalty_function} rather than one that directly addresses the consensus constraint  deterministically is given in Remark \ref{remark_penalty}, motivated by following the algorithm derivation. For future reference, we also define the local penalty as 
\begin{align}\label{eq:local_penalty_function}
\psi_{i,c}(f_i)\! &= \mbE_{\bbx_i, \bby_i}\! \Big[ \ell_i(f_i\big(\bbx_i), y_i\big)\!\Big]\! \!+\!\frac{\lambda}{2}\|f_i \|^2_{\ccalH} \! \nonumber \\
&\qquad\qquad\! \! + \! \frac{c}{2}  \!\!\sum_{j\in n_i}\!  \mathbb{E}_{\bbx_i}\left\{ [f_i(\bbx_i) \! - \! f_j(\bbx_i)]^2 \right\}
\end{align}
and we observe from \eqref{eq:penalty_function} - \eqref{eq:local_penalty_function} that $\psi_c(f) = \sum_i \psi_{i,c}(f_i)$. Further define $f^*_c =\argmin_{f\in\ccalH^V} \psi_c(f)$.
We note that in the vector-valued decision variable case, other techniques to address the constraint in \eqref{eq:main_prob} are possible such as primal-dual methods \cite{Koppel2015a} or dual methods \cite{Suzuki2013}, but the Representer Theorem has not been established for RKHS-valued stochastic saddle point problems. It is an open question whether expressions of the form \eqref{eq:kernel_expansion} apply to problems with general functional constraints, but this matter is beyond the scope of this work. Therefore, these other approaches which make use of Lagrange duality do not readily extend to the nonparametric setting considered here. 
\subsection{Functional Stochastic Gradient Method}\label{subsec:fsg}
%
\begin{algorithm}[t]
\caption{Greedy Projected Penalty Method}
\begin{algorithmic}
\label{alg:soldd}
\REQUIRE $\{\bbx_t,\bby_t,\eta_t,\epsilon_t \}_{t=0,1,2,...}$
\STATE \textbf{initialize} ${f}_{i,0}(\cdot) = 0, \bbD_{i,0} = [], \bbw_0 = []$, i.e. initial dictionary, coefficients are empty for each $i\in\ccalV$
\FOR{$t=0,1,2,\ldots$}
	\STATE {\bf loop in parallel} for agent $ i \in \ccalV$ 
	\INDSTATE Observe local training example realization $(\bbx_{i,t}, y_{i,t})$
	\INDSTATE{Send obs. $\bbx_{i,t}$ to nodes $j\in n_i$, receive scalar $f_{j,t}(\bbx_{i,t})$ }
	\INDSTATE{Receive obs. $\bbx_{j,t}$ from nodes $j \in n_i$, send $f_{i,t}(\bbx_{j,t})$ } 
	\INDSTATE Compute unconstrained stochastic grad. step [cf. \eqref{eq:sgd_tilde}]
	$$\tilde{f}_{i,t+1}(\cdot) = (1-\eta_t\lambda){f}_{i,t} - \eta_t\nabla_{f_i} \hat{\psi}_{i,c}(f_i(\bbx_{i,t}), \bby_{i,t})\; .$$\vspace{-4mm}
	\INDSTATE Update params: $\tbD_{i,t+1}\! =\! [\bbD_{i,t},\;\;\bbx_{i,t}]$, $\tbw_{i,t+1} $ [cf. \eqref{eq:param_tilde}]
		\INDSTATE Greedily compress function using matching pursuit \vspace{-1mm}
	$$\!\!\!\!\!({f}_{i,t+1}\!,\bbD_{i,t+1}\!,\!\bbw_{i,t+1}\!) \!= \textbf{KOMP}(\!\tilde{f}_{i,t+1},\!\tbD_{i,t+1},\!\tbw_{i,t+1},\!\epsilon_t\!)$$ 
	\STATE {\bf end loop}
\ENDFOR
\end{algorithmic}
\end{algorithm}

Given that the data distribution $\mathbb{P}(\bbx, \bby)$ is unknown, minimizing $\psi_c(f)$ directly via variational inference is not possible. Rather than postulate a specific distribution for $(\bbx, \bby)$, we only assume access to sequentially available (streaming) independent and identically distributed samples $(\bbx_t, \bby_t)$ from their joint density. Then, we may wield tools from stochastic approximation to minimize \eqref{eq:penalty_function}, which in turn yields a solution to \eqref{eq:main_prob}. 
Begin by defining, $\hat{\psi}_c(f(\bbx_t), \bby_t)$, the stochastic approximation of the penalty function $\psi_c(f)$, evaluated at a realization $(\bbx_t, \bby_t)$ of the stacked random pair $(\bbx, \bby)$:
\begin{align}\label{eq:penalty_function_hat}
\hat{\psi}_c(f(\bbx_t), \bby_t)\! &=  \!\sum_{i\in\ccalV}\!\Big( \ell_i(f_i\big(\bbx_{i,t}), y_{i,t}\big) \!+\!\frac{\lambda}{2}\|f_i \|^2_{\ccalH} \! \nonumber \\
&\qquad+ \! \frac{c}{2}  \!\sum_{j\in n_j} \!( f_i(\bbx_{i,t}) \! - \! f_j(\bbx_{i,t}))^2 \Big)
\end{align}
and the local instantaneous penalty function $\hat{\psi}_{i,c}(f_i(\bbx_{i,t}), \bby_{i,t})$ similarly. To compute the functional stochastic gradient of $\psi_c(f)$ evaluated at a sample point $(\bbx_t, \bby_t)$, we first address the local loss $ \ell_i(f_i\big(\bbx_{i,t}), y_{i,t})$ in \eqref{eq:penalty_function_hat} as \cite{Kivinen2004,POLK}:
\begin{align}\label{eq:stochastic_grad}
\! \! \!\nabla_{f_i} \ell_i(f_i(\bbx_{i,t}),y_{i,t})(\cdot)
= \frac{\partial \ell_i(f_i(\bbx_{i,t}),y_{i,t})}{\partial f_i(\bbx_{i,t})}\frac{\partial f_i(\bbx_{i,t})}{\partial f_i}(\cdot)
\end{align}
where we have applied the chain rule. Now, define the short-hand notation $$\ell_i'(f_i(\bbx_{i,t}),y_{i,t}): ={\partial \ell_i(f_i(\bbx_{i,t}),y_{i,t})}/{\partial f_i(\bbx_{i,t})} $$ for the derivative of $\ell_i(f(\bbx_{i,t}),y_{i,t})$ with respect to its first scalar argument $f_i(\bbx_{i,t})$ evaluated at $\bbx_{i,t}$. To evaluate the second term on the right-hand side of \eqref{eq:stochastic_grad}, differentiate both sides of the expression defining the reproducing property of the kernel [cf. \eqref{eq:rkhs_properties}(i)] with respect to $f_i$ to obtain
\begin{align}\label{eq:stochastic_grad2}
\frac{\partial  f_i(\bbx_{i,t})}{\partial f_i} = \frac{\partial \langle f_i , \kappa(\bbx_{i,t}, \cdot) \rangle _{\ccalH}}{\partial f_i}
= \kappa(\bbx_{i,t},\cdot)
\end{align}
Then, given \eqref{eq:stochastic_grad} - \eqref{eq:stochastic_grad2}, we may compute the overall gradient of the instantaneous penalty function $\hat{\psi}_c(f(\bbx_t), \bby_t)\!$ in \eqref{eq:penalty_function_hat} as
\begin{align}\label{eq:penalty_function_stoch_grad}
\nabla_f \hat{\psi}_c(f(\bbx_t), \bby_t)\! &= \text{vec}\Big[\ell_i'(f_i(\bbx_{i,t}),y_{i,t}) \kappa(\bbx_{i,t},\cdot)\!+\lambda f_i  \\
&\qquad+ \! c  \!\sum_{j\in n_i} ( f_i(\bbx_{i,t}) \! - \! f_j(\bbx_{i,t}))\kappa(\bbx_{i,t},\cdot) \Big] \nonumber
\end{align}
where on the right-hand side of \eqref{eq:penalty_function_stoch_grad}, we have defined the vector stacking notation $\text{vec}[\cdot]$ to denote the stacking of $V$ component-wise functional gradients, each associated with function $f_i$, $i \in \ccalV$, and used the fact that the variation of the instantaneous approximate of the cross-node term, $[f_i(\bbx_i) \! - \! f_j(\bbx_i)]^2$, by the same reasoning as \eqref{eq:stochastic_grad} - \eqref{eq:stochastic_grad2}, is $2 [f_i(\bbx_{i,t}) \! - \! f_j(\bbx_{i,t})]\kappa(\bbx_{i,t}, \cdot)$.
With this computation in hand, we present the stochastic gradient method for the $\lambda$-regularized multi-agent expected risk minimization problem in \eqref{eq:main_prob} as
\begin{align}\label{eq:sgd_hilbert}
f_{t+1} &=(1-\eta_t \lambda ) f_{t} - \eta_t \text{vec}\Big[ \ell_i'(f_{i,t}(\bbx_{i,t}),y_{i,t}) \kappa(\bbx_{i,t},\cdot)\nonumber \\
&\quad + \! c  \!\sum_{j\in n_i} ( f_{i,t}(\bbx_{i,t}) \! - \! f_{j,t}(\bbx_{i,t}))\kappa(\bbx_{i,t},\cdot) \Big]  \; ,
\end{align}
where $\eta_t> 0$ is an algorithm step-size either chosen as diminishing with $\ccalO(1/t)$ or a small constant -- see Section \ref{sec:convergence}. We may glean from \eqref{eq:sgd_hilbert} that the update for the network-wide function $f_t$ decouples into ones for each agent $i\in\ccalV$, using the node-separability of the penalty $\psi_c(f) = \sum_i \psi_{i,c}(f_i)$, i.e.,
\begin{align}\label{eq:sgd_hilbert_local}
f_{i,t+1} &=(1-\eta_t \lambda ) f_{i,t} - \eta_t \Big[\ell_i'(f_{i,t}(\bbx_{i,t}),y_{i,t}) \kappa(\bbx_{i,t},\cdot)\nonumber \\
&\quad + \! c  \!\sum_{j\in n_i} ( f_{i,t}(\bbx_{i,t}) \! - \! f_{j,t}(\bbx_{i,t}))\kappa(\bbx_{i,t},\cdot) \Big]  \; .
\end{align}
We further require that, given $\lambda > 0$, the step-size satisfies $\eta_t < 1/\lambda$ and the global sequence is initialized as $f_0 = 0 \in \ccalH^V$. With this initialization, the Representer Theorem \eqref{eq:kernel_expansion} implies that, at time $t$, the function $f_{i,t}$ admits an expansion in terms of feature vectors $\bbx_{i,t}$ observed thus far as
\begin{align}\label{eq:kernel_expansion_t}
f_{i,t}(\bbx) 
= \sum_{n=1}^{t-1} w_{i,n} \kappa(\bbx_{i,n}, \bbx)
= \bbw_{i,t}^T\boldsymbol{\kappa}_{\bbX_{i,t}}(\bbx) \; .
\end{align}
On the right-hand side of \eqref{eq:kernel_expansion_t} we have introduced the notation $\bbX_{i,t} = [\bbx_{i,1}, \ldots, \bbx_{i,t-1}]\in \reals^{p\times (t-1)}$, $\boldsymbol{\kappa}_{\bbX_{i,t}}(\cdot) = [\kappa(\bbx_{i,1},\cdot),\ \ldots\ ,\kappa(\bbx_{i,t-1},\cdot)]^T$, and $\bbw_{i,t} = [w_{i,1} , \ \ldots \, w_{i,t-1}] \in \reals^{t-1}$. Moreover, observe that the kernel expansion in \eqref{eq:kernel_expansion_t}, taken together with the functional update \eqref{eq:sgd_hilbert}, yields the fact that performing the stochastic gradient method in $\ccalH^V$ amounts to the following $V$ parallel parametric updates on the kernel dictionaries $\bbX_i$ and coefficients $\bbw_i$:
\begin{align}\label{eq:param_update} 
\bbX_{i,t+1} &= [\bbX_{i,t}, \;\; \bbx_{i,t}] \; ,  \\
\! [\bbw_{i,t+1}]_u \! \!& = \! \!
\begin{cases}
    (1 - \eta_t \lambda) [\bbw_{i,t}]_u \quad \text{for } 0 \leq u \leq t-1 \\
   \! \! -\eta_t\!\Big(\!\!\ell_i'(f_{i,t}(\!\bbx_{i,t}\!),\!y_{i,t}\!) \!\!+\!  c  \!\sum_{j\in n_i}\! ( f_{i,t}(\bbx_{i,t}\!) \! -\!\! f_{j,t}(\bbx_{i,t}\!)\!)\!\! \Big), 
  \end{cases} \nonumber 
\end{align}
where the second case on the last line of \eqref{eq:param_update} is for $u=t$.
Observe that this update causes $\bbX_{i,t+1}$ to have one more column than $\bbX_{i,t}$. We define the \emph{model order} as number of data points $M_{i,t}$ in the dictionary of agent $i$ at time $t$ (the number of columns of $\bbX_t$). FSGD is such that $M_{i,t}=t-1$, and hence grows unbounded with iteration index $t$. Next we address this intractable memory growth such that we may execute stochastic descent through low-dimensional projections of the stochastic gradient, inspired by \cite{POLK}. First, we clarify the motivation for the choice of the penalty function \eqref{eq:penalty_function}.
\begin{remark}\label{remark_penalty}\normalfont
In principle, it is possible to address the RKHS-valued consensus constraint in \eqref{eq:main_prob} directly, through primal-only  stochastic methods, by introducing the penalty function
\begin{align}\label{eq:penalty_function_remark}
\tilde{\psi}_c(f)\! &=  \!\!\sum_{i\in\ccalV}\!\!\Big(\!\mbE_{\bbx_i, \bby_i}\! \Big[\! \ell_i(f_i\big(\bbx_i), y_i\big)\!\Big]\! \!+\!\frac{\lambda}{2}\|f_i \|^2_{\ccalH}  \!+ \! \frac{c}{2}  \!\!\sum_{j\in n_i}\! \! \|f_i \!- \! f_j \|^2_{\ccalH}\! \Big)
\end{align}
Observe, however, that FSGD applied to \eqref{eq:penalty_function_remark}, using comparable reasoning to that which leads to \eqref{eq:sgd_hilbert_local} from \eqref{eq:penalty_function}, yields
\begin{align}\label{eq:sgd_hilbert_remark}
f_{i,t+1} &=(1-\eta_t \lambda ) f_{i,t} - \eta_t \Big[\nabla_{f_i} \ell_i'(f_{i,t}(\bbx_{i,t}),y_{i,t}) \kappa(\bbx_{i,t},\cdot)\nonumber \\
&\quad + \! c  \!\sum_{j\in n_i} ( f_{i,t} \! - \! f_{j,t}) \Big]  \; .
\end{align}
Unfortunately, we cannot inductively define a parametric representation of \eqref{eq:sgd_hilbert_remark} for node $i$ in terms of its own kernel dictionaries and weights independently of the \emph{entire function} associated to node $j$, since the last term in \eqref{eq:sgd_hilbert_remark} lives directly in the Hilbert space. Thus, to implement \eqref{eq:sgd_hilbert_remark} each agent would need to store the entire kernel dictionary and weights of all its neighbors at each step, which is impractically costly. The use of \eqref{eq:penalty_function} rather than \eqref{eq:penalty_function_remark} is further justified that under a hypothesis regarding the mean transformation of the local data spaces, $\mathbb{E}_{\bbx_i}[\kappa(\bbx_i, \cdot)]$, consensus with respect to the Hilbert norm, in addition to the mean square sense, is achieved when the penalty coefficient is $c\rightarrow \infty$ (see Section \ref{sec:convergence} for details).
\end{remark}
\subsection{Sparse Subspace Projections}\label{subsec:proj}
To mitigate the complexity growth noted in Section \ref{subsec:fsg}, we approximate the function sequence \eqref{eq:sgd_hilbert} by one that is orthogonally projected onto subspaces $\ccalH_\bbD \subseteq \ccalH$ that consist only of functions that can be represented using some dictionary $\bbD = [\bbd_1,\ \ldots,\ \bbd_M] \in \reals^{p \times M}$, i.e., $\ccalH_\bbD = \{f\ :\ f(\cdot) = \sum_{n=1}^M w_n\kappa(\bbd_n,\cdot) = \bbw^T\boldsymbol{\kappa}_{\bbD}(\cdot) \}=\text{span}\{\kappa(\bbd_n, \cdot) \}_{n=1}^M$, and $\{\bbd_n\} \subset \{\bbx_u\}_{u\leq t}$. For convenience we define $[\boldsymbol{\kappa}_{\bbD}(\cdot)=\kappa(\bbd_1,\cdot) \ldots \kappa(\bbd_M,\cdot)]$, and $\bbK_{\bbD,\bbD}$ as the resulting kernel matrix from this dictionary. We enforce function parsimony by selecting dictionaries $\bbD_i$ with $M_{i,t} << \ccalO(t)$ for each $i$ \cite{POLK}.

To be specific, we propose replacing the local update \eqref{eq:sgd_hilbert_local} in which the dictionary grows at each iteration by its projection onto subspace $\ccalH_{\bbD_{i,t+1}}=\text{span}\{ \kappa(\bbd_{i,n}, \cdot) \}_{n=1}^{M_{t+1}}$ as
\begin{align}\label{eq:projection_hat}
{f}_{i,t+1}& = \!\!\argmin_{f \in \ccalH_{\bbD_{i,t+1}}}\! \!\! \Big\lVert f \!-\!\left(\!f_{i,t}  \!-\! \eta_t \nabla_{f_i}\hat{\psi}_{i,c}(f_i(\bbx_{i,t}), y_{i,t})\right) \! \!
\Big\rVert_{\ccalH}^2 \nonumber \\
&:=\ccalP_{\ccalH_{\bbD_{i,t+1}}} \Big[ 
(1-\eta_t \lambda) f_{i,t} 
- \eta_t\Big( \nabla_{f_i} \ell_i(f_{i,t}(\bbx_{i,t}),y_{i,t}) \nonumber \\
& \qquad\quad+ \! c  \!\sum_{j\in n_i} ( f_{i,t}(\bbx_{i,t}) \! - \! f_{j,t}(\bbx_{i,t}))\kappa(\bbx_{i,t},\cdot)\Big)\Big] . 
\end{align}
where we define the projection operator $\ccalP$ onto subspace $\ccalH_{\bbD_{i,t+1}}\subset \ccalH$ by the update \eqref{eq:projection_hat}.

{\bf Coefficient update} The update \eqref{eq:projection_hat}, for a fixed dictionary $\bbD_{i,t+1} \in \reals^{p\times M_{t+1}}$, yields one in the coefficient space only. This fact may be observed by defining the un-projected stochastic gradient step starting at function ${f}_{i,t}$ parameterized by dictionary $\bbD_{i,t}$ and coefficients $\bbw_{i,t}$:
\begin{align}\label{eq:sgd_tilde}
\tilde{f}_{i,t+1} ={f}_{i,t} -\eta_t \nabla_{f_i} \hat{\psi}_{i,c}(f_i(\bbx_{i,t}), y_{i,t}) \; .
\end{align}
This update may be represented using dictionary and weights
\begin{align}\label{eq:param_tilde} 
\tbD_{i,t+1} & = [\bbD_{i,t},\;\;\bbx_{i,t}]\; ,  \\
\! [\tbw_{i,t+1}]_u \! \!& = \! \!
\begin{cases}
    (1 - \eta_t \lambda) [\bbw_{i,t}]_u \quad \text{for } 0 \leq u \leq t -1\\
   \! \! -\eta_t\!\Big(\!\!\ell_i'(f_{i,t}(\!\bbx_{i,t}\!),\!y_{i,t}\!) \!\!+\!  c  \!\sum_{j\in n_i}\! ( f_{i,t}(\bbx_{i,t}\!) \! -\!\! f_{j,t}(\bbx_{i,t}\!)\!)\!\! \Big), 
  \end{cases} \nonumber 
\end{align}
where the last coefficient is for $u=t$.
%
Note that $\tbD_{i,t+1}$ has $\tilde{M}=M_{i,t} + 1$ columns, which is also the length of $\tbw_{i,t+1}$. For a fixed $\bbD_{i,t+1}$, the stochastic projection \eqref{eq:projection_hat} is a least-squares update on the coefficient vector:  the Representer Theorem allows us to rewrite \eqref{eq:projection_hat} in terms of kernel expansions as in Section 3.2 of \cite{POLK}, which yields
\begin{equation} \label{eq:hatparam_update}
\bbw_{i,t+1}=  \bbK_{\bbD_{i,t+1} \bbD_{i,t+1}}^{-1} \bbK_{\bbD_{i,t+1} \tbD_{i,t+1}} \tbw_{i,t+1} \;,
\end{equation}
where we define the cross-kernel matrix $\bbK_{\bbD_{i,t+1},\tbD_{i,t+1}}$ whose $(n,m)^\text{th}$ entry is given by $\kappa(\bbd_{i,n},\tbd_{i,m})$. The other kernel matrices $\bbK_{\tbD_{i,t+1},\tbD_{i,t+1}}$ and $\bbK_{\bbD_{i,t+1},\bbD_{i,t+1}}$ are defined similarly. Observe that ${M}_{i,t+1}$ is the number of columns in $\bbD_{i,t+1}$, while $\tilde{M}_i=M_{i,t} + 1$ is the number of columns in $\tbD_{t+1}$ [cf. \eqref{eq:param_tilde}]. Given that the local projections of $\tilde{f}_{i,t+1}$ onto stochastic subspaces $\ccalH_{\bbD_{i,t+1}}$, for a fixed node-specific dictionaries $\bbD_{i,t+1}$, is a least-squares problem, we now detail the kernel dictionary $\bbD_{i,t+1}$ selection from past data $\{\bbx_{i,u}, y_{i,u}\}_{u \leq t}$.

{\bf Dictionary Update} The selection procedure for the kernel dictionary $\bbD_{i,t+1}$ is based upon greedy compression \cite{needell2008greedy}: function $\tilde{f}_{i,t+1}$ defined by the stochastic gradient method without projection is parameterized by dictionary $\tbD_{i,t+1}$ [cf. \eqref{eq:param_tilde}] of model order $\tilde{M_i}={M}_{i,t} +1$. We form $\bbD_{i,t+1}$ by selecting a subset of $M_{i,t+1}$ columns from $\tbD_{i,t+1}$ that best approximate $\tilde{f}_{i,t+1}$ in terms of Hilbert norm error, which may be done by executing \emph{kernel orthogonal matching pursuit} (KOMP) \cite{Pati1993,Vincent2002} with error tolerance $\epsilon_t$ to find a kernel dictionary matrix $\bbD_{i,t+1}$ based on the one which adds the latest sample point $\tbD_{i,t+1}$. This choice is due to the fact that we can tune its stopping criterion to guarantee stochastic descent, and guarantee the model order of the learned function remains finite -- see Section \ref{sec:convergence} for details.

We now describe the variant of KOMP we propose using, called Destructive KOMP with Pre-Fitting (see \cite{Vincent2002}, Section 2.3). Begin with an input a candidate function $\tilde{f}$ of model order $\tilde{M}$ parameterized by kernel dictionary $\tbD\in\reals^{p\times\tilde{M}}$ and coefficients $\tbw\in\reals^{\tilde{M}}$. The method then approximates $\tilde{f}$ by a function $f\in \ccalH$ with a lower model order. Initially, this sparse approximation is the original function $f = \tilde{f} $ so that its dictionary is initialized with that of the original function $\bbD=\tbD$, with corresponding coefficients  $\bbw=\tbw$. 
Then, the algorithm sequentially removes dictionary elements from the initial dictionary $\tbD$, yielding a sparse approximation $f$ of $\tilde{f}$, until the error threshold $\|f - \tilde{f} \|_{\ccalH} \leq \eps_t $ is violated, in which case it terminates. See Appendix A for further details.


We summarize the key steps of the proposed method in Algorithm \ref{alg:soldd} for solving \eqref{eq:main_prob} while maintaining a finite model order, thus allowing for the memory-efficient learning of nonparametric regression functions online in multi-agent systems. The method, Greedy Projected Penalty Method, executes the stochastic projection of the functional stochastic gradient iterates onto sparse subspaces $\ccalH_{\bbD_{i,t+1}}$ stated in \eqref{eq:projection_hat}. Initial functions are set to null $f_{i,0}=0$, i.e., it has empty dictionary $\bbD_{i,0}=[]$ and coefficient vector $\bbw_{i,0}=[]$. The notation $[]$ is used to denote the empty matrix or vector respective size $p\times0$ or $0$. Then, at each step, given an independent training example $(\bbx_{i,t}, y_{i,t})$ and step-size $\eta_t$, we compute the \emph{unconstrained} functional stochastic gradient iterate \eqref{eq:sgd_tilde} with respect to the instantaneous penalty function \eqref{eq:penalty_function_hat} which admits the parameterization $\tbD_{i,t+1}$ and $\tbw_{i,t+1}$ as stated in \eqref{eq:param_tilde}. These parameters are then fed into KOMP with approximation budget $\eps_t$, such that $(f_{i,t+1}, \bbD_{i,t+1}, \bbw_{i,t+1})= \text{KOMP}(\tilde{f}_{i,t+1},\tilde{\bbD}_{i,t+1}, \tilde{\bbw}_{i,t+1},\eps_t)$.

%


\section{Convergence Analysis}\label{sec:convergence}

We turn to establishing that the method presented in Algorithm \ref{alg:soldd} converges with probability $1$ to the minimizer of the penalty function $\psi_c(f)$ [cf. \eqref{eq:penalty_function}] when attenuating algorithm step-sizes are used, and to a neighborhood of the minimizer along a subsequence when constant step-sizes are used. Moreover, for the later case, the kernel dictionary that parameterizes the regression function $f_i$ for each agent $i$ remains finite in the worst case. This analysis is an application of Section IV of \cite{POLK}, but these results, together with the properties of the penalty function $\psi_c(f)$ allow us to establish bounds on the deviation for each individual in the network from the common globally optimal regression function.

%
%

Before analyzing the proposed method developed in Section \ref{sec:algorithm}, we define key quantities to simplify the analysis and introduce standard assumptions which are necessary to establish convergence. 
Define the local projected stochastic functional gradient associated with the update in \eqref{eq:projection_hat} as
\begin{align}\label{eq:proj_grad}
\tilde{\nabla}_{f_i} \!\hat{\psi}_{i,c} &(f_{i,t}(\bbx_{i,t}\!),y_{i,t})=\! \\
& \!\!\Big(\! f_{i,t} - \ccalP_{ \ccalH_{\bbD_{i,t+1}}}\! \!\Big[\!
 f_{i,t} 
- {\eta}_t \!{\nabla}_{f_i}\hat{\psi}_{i,c}(f_{i,t}(\bbx_{i,t}\!),\!y_{i,t}\!) \!\Big]\!\Big)\!/\eta_t \nonumber
\end{align}
such that the local update of Algorithm \ref{alg:soldd} [cf. \eqref{eq:projection_hat}] may be expressed as a stochastic descent using projected functional gradients
%
$f_{i,t+1} = f_{i,t} - \eta_t \tilde{\nabla}_{f_i} \!\hat{\psi}_{i,c} (f_{i,t}(\bbx_{i,t}\!),y_{i,t})\; .$
%
The definitions of \eqref{eq:proj_grad} and the local stochastic gradient ${\nabla}_{f_i} \!\hat{\psi}_{i,c} (f_{i,t}(\bbx_{i,t}\!),y_{i,t})$ may be stacked to analyze the global convergence behavior of the algorithm. For further reference, we define the stacked projected functional stochastic gradient of the penalty function as $\tilde{\nabla}_{f} \hat{\psi}_{c} (f_{t}(\bbx_{t}),\bby_{t}) =[ \tilde{\nabla}_{f_1} \hat{\psi}_{1,c} (f_{1,t}(\bbx_{1,t}\!),y_{1,t}) ; \cdots ; \tilde{\nabla}_{f_V} \hat{\psi}_{V,c} (f_{V,t}(\bbx_{V,t}\!),y_{V,t})]$. Then the stacked global update of the algorithm is
\begin{equation}\label{eq:iterate_tilde}
f_{t+1} = f_{t} - \eta_t \tilde{\nabla}_{f} \hat{\psi}_{c} (f_{t}(\bbx_{t}),\bby_{t})\; .
\end{equation}
Moreover, observe that the stochastic functional gradient in \eqref{eq:penalty_function_stoch_grad}, based upon the fact that $(\bbx_t, y_t)$ are independent and identically distributed realizations of the random pair $(\bbx, y)$, is an unbiased estimator of the true functional gradient of the penalty function $\psi_c(f)$ in \eqref{eq:penalty_function}, i.e.
\begin{equation}\label{eq:unbiased}
\mathbb{E}[\nabla_f \hat{\psi}_c(f(\bbx_t), \bby_t) \given \ccalF_t ] =\nabla_f {\psi}_c(f)
\end{equation}
for all $t$. In \eqref{eq:unbiased}, we denote as $\ccalF_t$ the sigma algebra which measures the algorithm history for times $u\leq t$, i.e. $\ccalF_t=\{\bbx_u, y_u, u_u\}_{u=1}^t$.
Next, we formally state technical conditions on the loss functions, data domain, and stochastic approximation errors that are necessary to establish convergence.
%

\begin{assumption}\label{as:first}
The feature space $\ccalX\subset\reals^p$ and target domain $\ccalY\subset\reals$ are compact, and the  kernel map may be bounded as
\begin{equation}\label{eq:bounded_kernel}
\sup_{\bbx\in\ccalX} \sqrt{\kappa(\bbx, \bbx )} = X < \infty
\end{equation}
%
%
\end{assumption}
%
%
\begin{assumption}\label{as:3}
The local losses $\ell_i(f_i(\bbx), y)$ are convex and differentiable with respect to the first (scalar) argument $f_i(\bbx)$ on $\reals$ for all $\bbx\in\ccalX$ and $y\in\ccalY$.  
Moreover, the instantaneous losses $\ell_i: \ccalH \times \ccalX \times \ccalY \rightarrow \reals$ are $C_i$-Lipschitz continuous for all $z \in \reals $ for a fixed $y\in\ccalY$
\begin{equation}\label{eq:lipschitz}
| \ell_i(z, y) - \ell_i( z', y ) | \leq C_i |z - z'|
\end{equation}
with $C :=\max_i C_i$ as the largest modulus of continuity.

\end{assumption}
%
%
\begin{assumption}\label{as:last}
The projected functional gradient of the instantaneous penalty function defined by stacking \eqref{eq:proj_grad} has finite conditional second moments:
\begin{equation}\label{eq:stochastic_grad_var}
\mathbb{E} [ \| \tilde{\nabla}_{f} \hat{\psi}_{c} (f_{t}(\bbx_{t}),\bby_{t})\|^2_{\ccalH} \mid \ccalF_t ] \leq \sigma^2
\end{equation}
\end{assumption}
%
Assumption \ref{as:first} holds in most settings by the data domain itself, and justifies the bounding of the loss. Taken together, these conditions permit bounding the optimal function $f^*_c$ in the Hilbert norm, and imply that the worst-case model order is guaranteed to be finite. Variants of Assumption \ref{as:3} appear in the analysis of stochastic descent methods in the kernelized setting \cite{Pontil05erroranalysis,1715525}, and is satisfied for supervised learning problems such as logistic regression, support vector machines with the square-hinge-loss, the square loss, among others. Moreover, it is standard in the analysis of descent methods (see \cite{nesterov1998introductory}). Assumption \ref{as:last} is common in stochastic methods, and ensures that the stochastic approximation error has finite variance. 

Next we establish a few auxiliary results needed in the proof of the main results. 
%
Specifically, we introduce a proposition which quantifies the error due to sparse projections in terms of the ratio of the compression budget to the learning rate.
%
\begin{proposition}\label{prop_projection}
Given independent realizations $(\bbx_t, \bby_t)$ of the random pair $(\bbx, \bby)$, the difference between the stacked projected stochastic functional gradient and the its un-projected variant defined by \eqref{eq:proj_grad} and \eqref{eq:penalty_function_stoch_grad}, respectively, is bounded as
\begin{equation}\label{eq:prop_projection}
 \| \tilde{\nabla}_{f} \hat{\psi}_{c} (f_{t}(\bbx_{t}),\bby_{t}) - \nabla_f \hat{\psi}_c(f(\bbx_t), \bby_t) \|_{\ccalH} \leq \frac{\eps_t V}{\eta_t}
\end{equation}
where $\eta_t>0$ denotes the algorithm step-size and $\eps_t>0$ is the approximation budget parameter of Algorithm \ref{alg:komp}.
\end{proposition}
\begin{myproof} See Appendix B. \end{myproof}

With the error induced by sparse projections quantified, we may now shift focus to analyzing the Hilbert-norm sub-optimality of the stacked iterates generated by Algorithm \ref{alg:soldd}. Specifically, we have a descent property of the sequence $\{f_t\}$.
\begin{lemma}(Stochastic Descent) \label{lemma1}
Consider the sequence generated $\{f_t \}$ by Algorithm \ref{alg:soldd} with $f_0 = 0$. Under Assumptions \ref{as:first}-\ref{as:last}, the following expected descent relation holds.
\begin{align}\label{eq:lemma1}
\E{\| f_{t+1} - f^*_c \|_{\ccalH}^2 \given \ccalF_t}
	 &\leq\| f_t \!- \!f^*_c \|_{\ccalH}^2
	\! -\! 2 \eta_t [\psi_c(f_t) \!-\!\psi_c(f^*_c) ]  \nonumber \\
	&\qquad + 2 \eps_t V \| f_t - f^*_c \|_{\ccalH}\! +\! \eta_t^2 \sigma^2
\end{align}
\end{lemma}
\begin{myproof} See Appendix B. \end{myproof}

Now that Lemma \ref{lemma1} establishes a descent-like property, we may apply the proof of Theorem 1 in \cite{POLK} to $\|f_t - f_c^*\|_{\ccalH}$ with diminishing step-sizes. Thus we have the following corollary.
\begin{corollary}\label{corollary:convergence}
Consider the sequence $\{f_t \}$ generated by Algorithm \ref{alg:soldd} with $f_0 = 0$ and regularizer $\lambda>0$. Under Assumptions \ref{as:first}-\ref{as:last} and the hypothesis that the projection sets $\ccalH_{\bbD_{i,t}}$ in \eqref{eq:projection_hat} are intersected with some finite Hilbert-norm ball $\|f\|_{\ccalH} \leq D$ for all $t$, with diminishing step-sizes and compression budget,  i.e.,
\begin{equation}\label{corollary_diminishing_stepsize}
\sum_{t=0}^\infty \eta_t = \infty \; , \quad \sum_{t=0}^\infty \eta_t^2 < \infty \; ,\quad \eps_t = \eta_t^2 \;, 
\end{equation}
such that $\eta_t<1/\lambda$, the sequence converges exactly to the minimizer of the penalty [cf. \eqref{eq:penalty_function}]: $f_t \rightarrow f_c^*$ with probability $1$.
\end{corollary}

To attain exact convergence to the minimizer of the penalty, $f^*_c$, we require the compression budget determining the error $\eps_t$ incurred by sparse projections to approach null. This means that to have exact convergence, we require the function representation to require an increasing amount of memory which is, in the limit, of infinite complexity. In contrast, when constant step-size and compression budget are used, then the algorithm settles to a neighborhood, as we state next.
\begin{thm}\label{corollary:convergence_constant}
The sequence $\{f_t \}$ generated by Algorithm \ref{alg:soldd} with $f_0 = 0$ and regularizer $\lambda>0$, under  Assumptions \ref{as:first}-\ref{as:last}, with constant step-size selection $\eta_t=\eta<1/\lambda$ and constant compression budget $\eps_t = \eps=K \eta^{3/2}$ for a positive constant $K$, converges to a neighborhood of $f_c^*$ with probability $1$:
\begin{equation}\label{corollary_constant_stepsize}
\!\!\! \!\!\!\!\! \liminf_t \! \|f_t \!-\! f_c^* \|_{\ccalH} \! \leq \!\! \frac{\sqrt{\eta}}{\lambda}\! \! \left[\!KV\!\! +\!\! \sqrt{\!K^2V^2 \!+\! \lambda \sigma^2}\right]\!\! = \! \ccalO(\!\sqrt{\eta}) \text{ a.s. } \!\!\!
\end{equation}
\end{thm}

\begin{myproof} See Appendix D. \end{myproof}

Empirically, the use of constant step-sizes has the effect of maintaining consistent algorithm adaptivity in the face of new data, at the cost of losing exact convergence. But this drawback is more than compensated for by the fact that in this case we may apply Theorem 3 of \cite{POLK}, which guarantees the model order of the function sequence remains finite, and in the worst case, is related to the covering number of the data domain

  %
\begin{figure*}%
\centering
\hspace{-8mm}
\begin{subfigure}{0.34\columnwidth}
\includegraphics[width=1\linewidth, height = 0.55\linewidth]{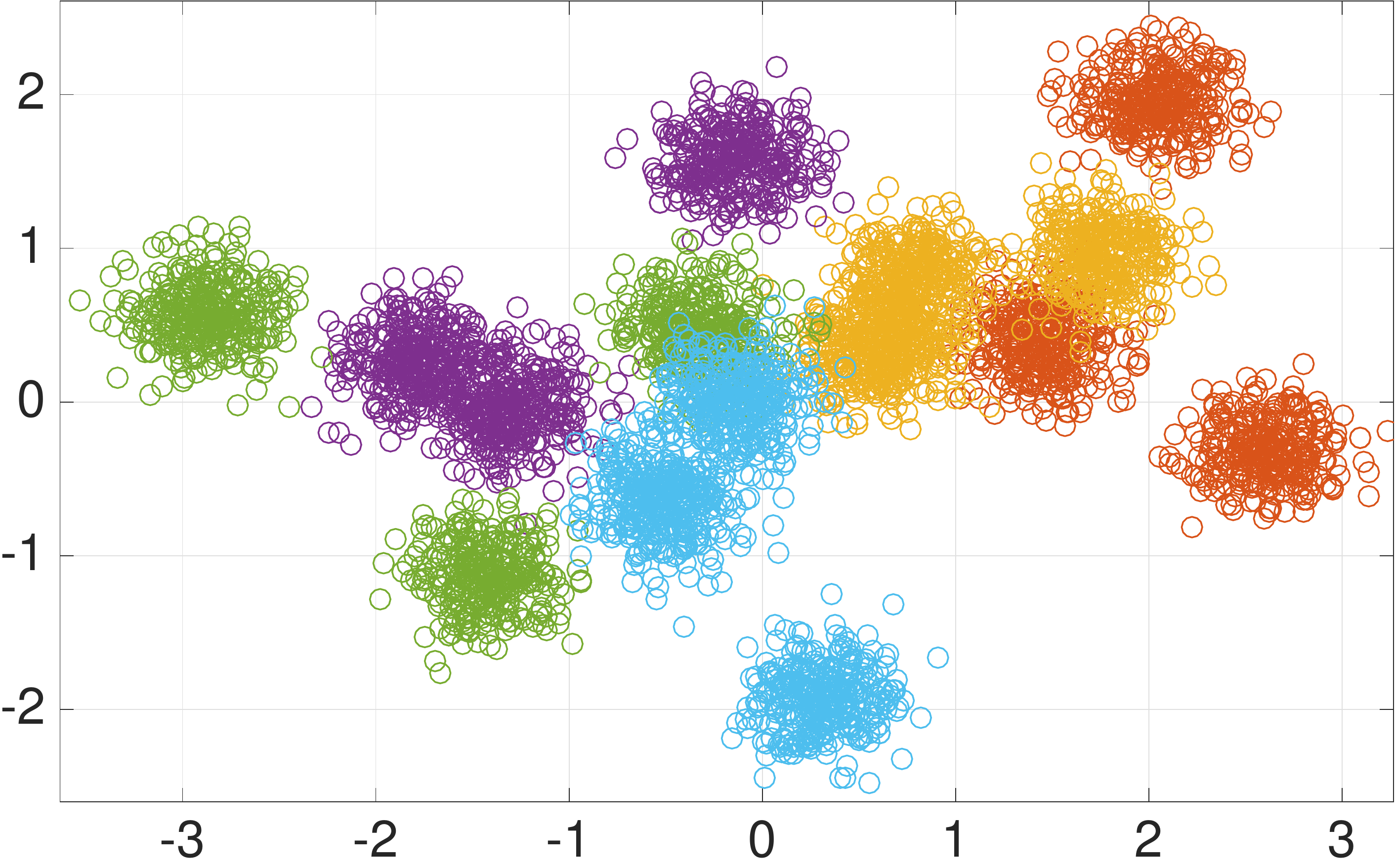}%
\caption{Gaussian Mixtures data.}%
\label{ch7_subfiga_gmm}%
\end{subfigure}
\begin{subfigure}{0.34\columnwidth}
\includegraphics[width=1.1\linewidth,height = 0.63\linewidth]{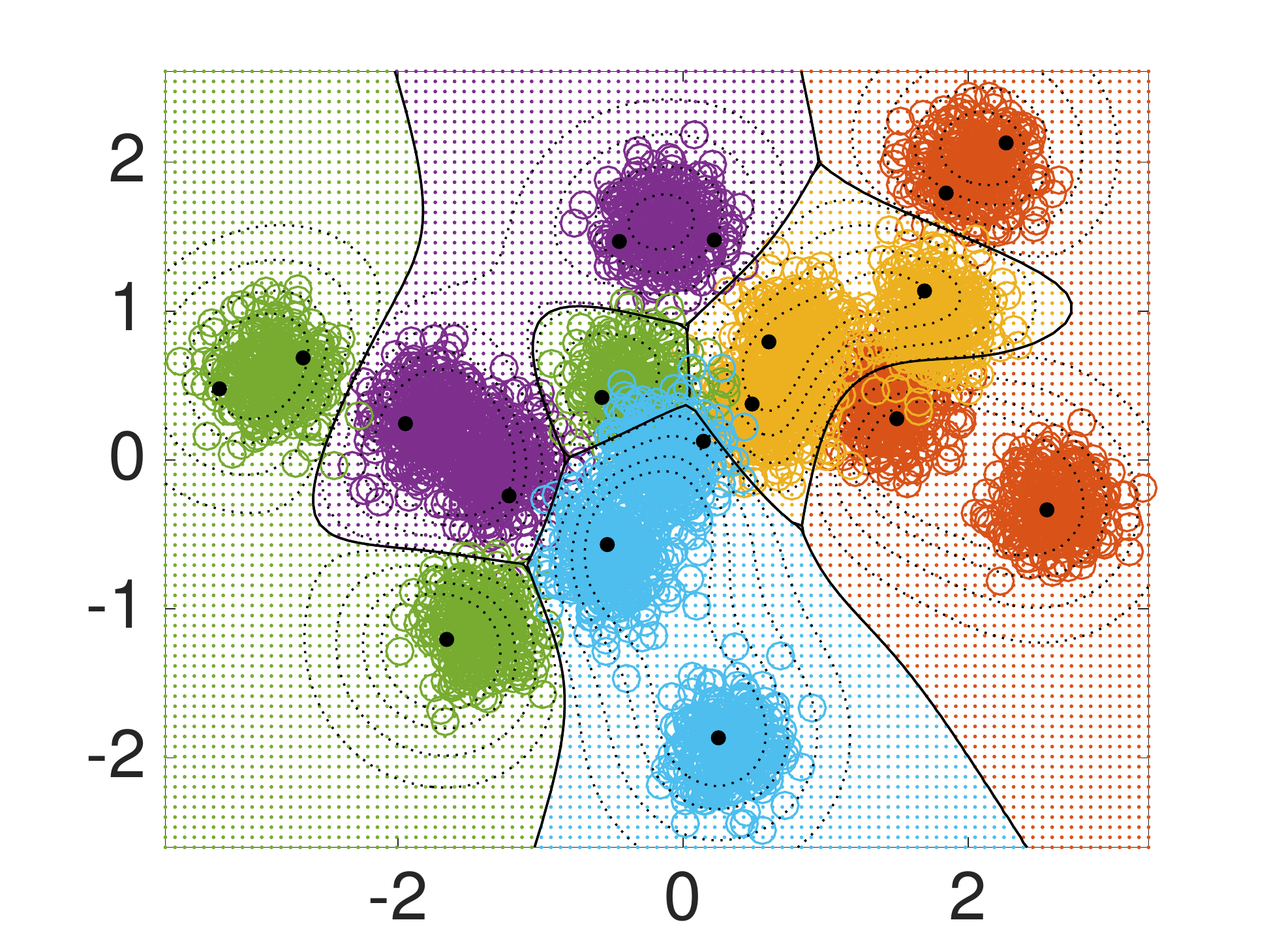}%
\caption{Logistic Decision surface.}%
\label{ch7_subfigb_gmm}%
\end{subfigure}%
\begin{subfigure}{0.34\columnwidth}
\includegraphics[width=1.1\linewidth,height = 0.63\linewidth]{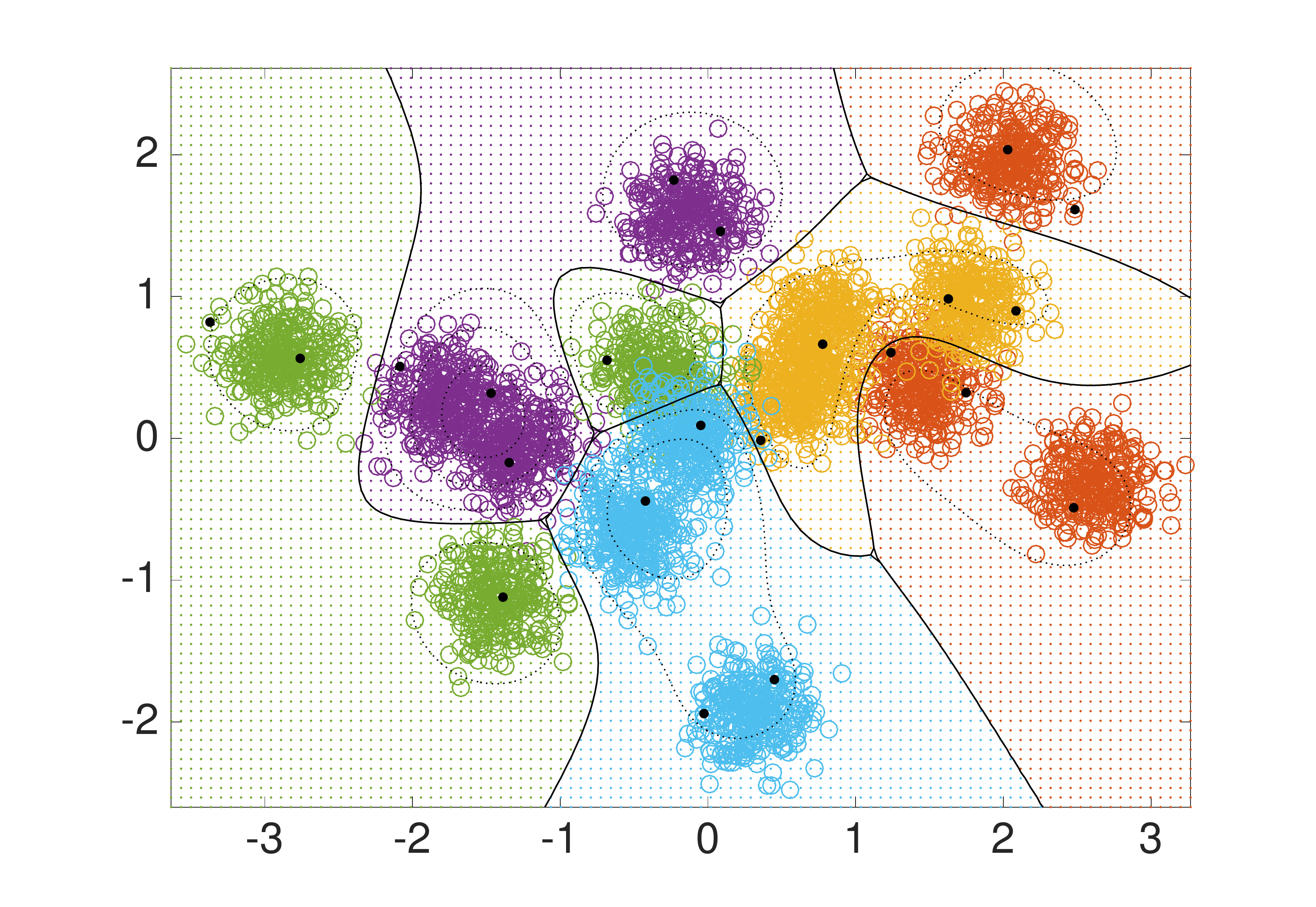}%
\caption{Hinge Decision surface.}%
\label{ch7_subfigc_gmm}%
\end{subfigure}%
\caption[Visualizations of decentralized adaptive accurate statistical learning.]{Visualizations of the Gaussian mixture data set (Figure \ref{ch7_subfiga_gmm}) as in \cite{Zhu2005} and the learned low-memory multi-class kernel logistic regressor of a randomly chosen agent in the network (Figure \ref{ch7_subfigb_gmm}), which attains $95.2\%$ classification accuracy on a hold-out test set. Curved black lines denote decision boundaries between classes; dotted lines denote confidence intervals; bold black dots denote kernel dictionary elements associated to an arbitrary $i\in\ccalV$. Kernel dictionary elements concentrate at peaks of the Gaussian clusters and near points of overlap between classes. In Figure \ref{ch7_subfigc_gmm} we plot the resulting decision surface learned by kernel SVM which attains $95.7\%$ accuracy -- the state of the art.}\label{ch7_fig:gmm}
\end{figure*}
\begin{corollary}\label{corollary:model_order}
Denote $f_t \in \ccalH^V$ as the stacked function sequence defined by Algorithm \ref{alg:soldd} with constant step-size $\eta_t=\eta<1/\lambda$ and approximation budget $\eps=K \eta^{3/2}$ where $K>0$ is an arbitrary positive scalar. Let $M_t$ be the model order of the stacked function $f_t$ i.e., the number of columns of the dictionary $\bbD_t$ which parameterizes $f_t$. Then there exists a finite upper bound $M^\infty$ such that, for all $t\geq0$, the model order is always bounded as $M_t\leq M^\infty$. 
\end{corollary}

Thus, only constant step-sizes attain a reasonable tradeoff between performance relative to $f_c^*$ and the complexity of storing the function sequence $\{f_t\}$: in this setting, we obtain approximate convergence to $f_c^*$ while ensuring the memory requirements are always finite, as stated in Corollary \ref{corollary:model_order}.

We are left to analyze the goodness of the solution $f^*_c$ as an approximation of the solution of the original problem \eqref{eq:main_prob}. In particular, we establish consensus in the mean square sense. Let us start by establishing that the penalty term is bounded by a $p^*/c$, where $p^*$ is the primal value of the optimization problem \eqref{eq:main_prob} and $c$ is the barrier parameter introduced in \eqref{eq:penalty_function}. 
%
\begin{proposition}\label{prop_constraint_violation}
  Let Assumptions  \ref{as:first} - \ref{as:last}  hold. Let $f^*_c$ be the minimizer of the penalty function \eqref{eq:penalty_function} and let $p^*$ be the primal optimal value of \eqref{eq:main_prob}. Then, it holds that 
  \begin{align}\label{eqn_constraint_bound}
\! \frac{1}{2} \sum_{i\in \ccalV} \!\!\sum_{j\in n_i}\!  \mathbb{E}_{\bbx_i}\left\{ [f^*_{c,i}(\bbx_i) \! - \! f^*_{c,j}(\bbx_i)]^2 \right\} \leq\frac{p^*}{c}.
    \end{align}
\end{proposition}
\begin{myproof}
See Appendix E.
  \end{myproof}

Proposition \ref{prop_constraint_violation} establishes a relationship between the choice of penalty parameter $c$ and constraint satisfaction. This result may be used to attain convergence in mean square of each individual agent's regression function to ones which coincide with one another. Under an additional hypothesis, we obtain exact consensus, as we state next.
\begin{thm}\label{thm_constraint_violation}
  Let Assumptions \ref{as:first} - \ref{as:last} hold. Let $f^*_c$ be the minimizer of the penalty function \eqref{eq:penalty_function}. Then, suppose the penalty parameter $c$ in \eqref{eq:penalty_function} approaches infinity $c\rightarrow \infty$, and that the node-pair differences $f_{i,c}^* - f_{j,c}^*$ are not orthogonal to mean transformation $\mathbb{E}_{\bbx_i}[\kappa(\bbx_i, \cdot)]$ of the local input spaces $\bbx_i$ for all $(i,j)\in\ccalE$. Then $f_{i,c}^* = f_{j,c}^*$ for all $(i,j)\in\ccalE$.
\end{thm}

\begin{myproof}
As a consequence, the limit of \eqref{eqn_constraint_bound} when $c$ tends to infinity yields consensus in $L^2$ sense, i.e.,
\begin{align}\label{eqn_consensus_1}
\lim_{c\to\infty} \, \frac{1}{2} \sum_{i\in \ccalV} \sum_{j\in n_i}\!  \mathbb{E}_{\bbx_i}\left\{ [f^*_{c,i}(\bbx_i) \! - \! f^*_{c,j}(\bbx_i)]^2 \right\} =0,
    \end{align}
which, by pulling the limit outside the sum in \eqref{eqn_consensus_1}, yields
\begin{align}\label{eqn_consensus}
\lim_{c\to\infty} \, \mathbb{E}_{\bbx_i}\left\{ [f^*_{c,i}(\bbx_i) \! - \! f^*_{c,j}(\bbx_i)]^2 \right\} =0 \; ,
    \end{align}
for all $(i,j)\in\ccalE$. Consensus in the mean square sense is a less stringent constraint that equality in the Hilbert norm as desired in \eqref{eq:main_prob}. In particular, for any $(i,j)\in\ccalE$, if $f_i=f_j,$ then consensus in the mean square sense is satisfied as well. Then, apply the reproducing property of the kernel \eqref{eq:rkhs_properties}(i), to write
  \begin{align}\label{eq:hilbert_norm_consensus1}
    0= \lim_{c\to\infty} \,  \mathbb{E}_{\bbx_i}\left\{ \left|<f^*_{c,i} \! - \! f^*_{c,j},k(\bbx_i,\cdot)> \right|\right\} \\
    \geq \lim_{c\to\infty} \,  \left| \mathbb{E}_{\bbx_i}\left\{ <f^*_{c,i} \! - \! f^*_{c,j},k(\bbx_i,\cdot)> \right\}\right| \nonumber \\
    =  \lim_{c\to\infty} \,  \left|  <f^*_{c,i} \! - \! f^*_{c,j},\mathbb{E}_{\bbx_i}k(\bbx_i,\cdot) >\right|  \nonumber 
  \end{align}
  where in the previous expression we pull the absolute value outside the expectation, and in the later we apply linearity of the expectation. Thus, \eqref{eq:hilbert_norm_consensus1} implies consensus is achieved with respect to the Hilbert norm, whenever the function differences $f^*_{c,i} \! - \! f^*_{c,j}$ are not orthogonal to $\mathbb{E}_{\bbx_i}[\kappa(\bbx_i, \cdot)]$, the mean of the transformation of the local input data $\bbx_i$. 
\end{myproof}

\begin{figure*}%
\centering 
\begin{subfigure}{0.33\columnwidth}
\includegraphics[width=1.05\linewidth, height = 0.73\linewidth]{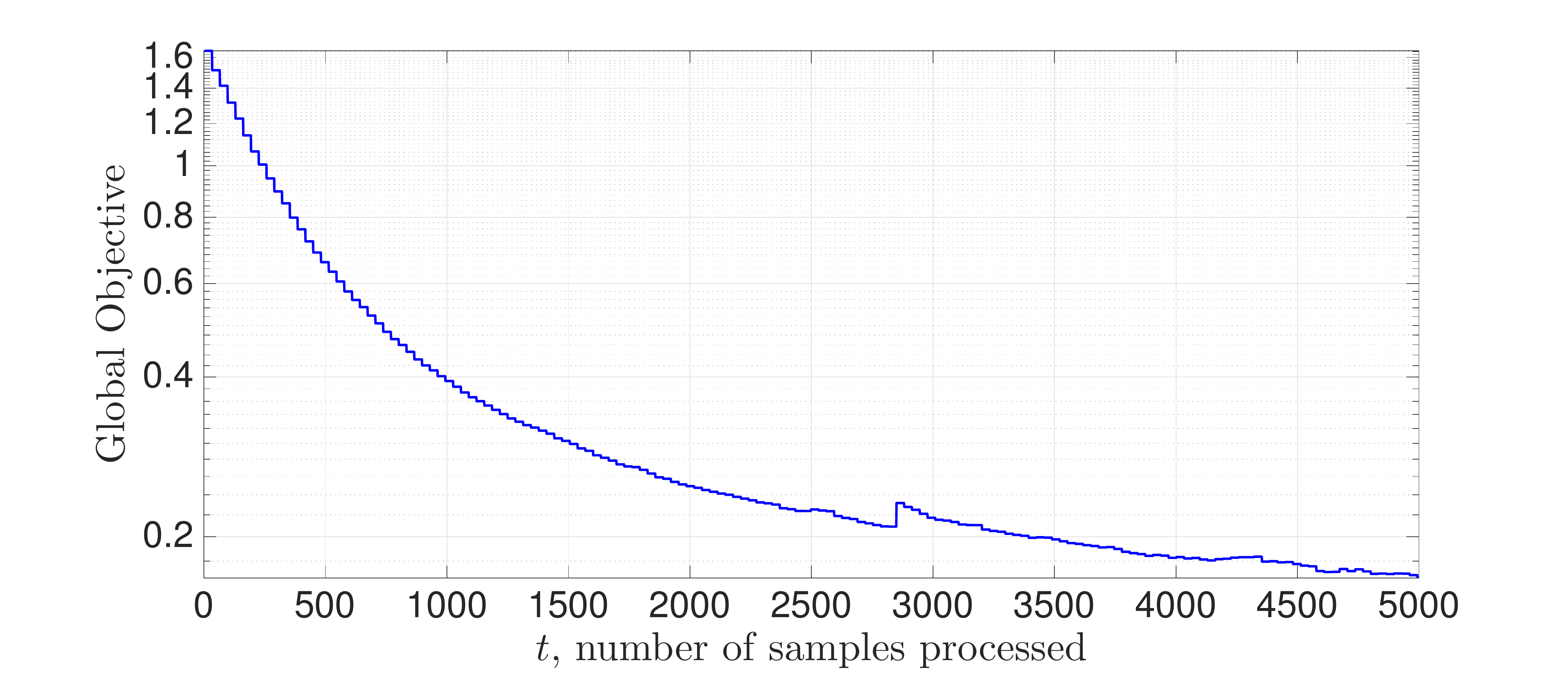}%
\caption{Global objective vs. samples processed}%
\label{subfiga_sample_path}%
\end{subfigure}\hspace{.5mm}
\begin{subfigure}{0.32\columnwidth}
\includegraphics[width=1.05\linewidth,height = 0.73\linewidth]{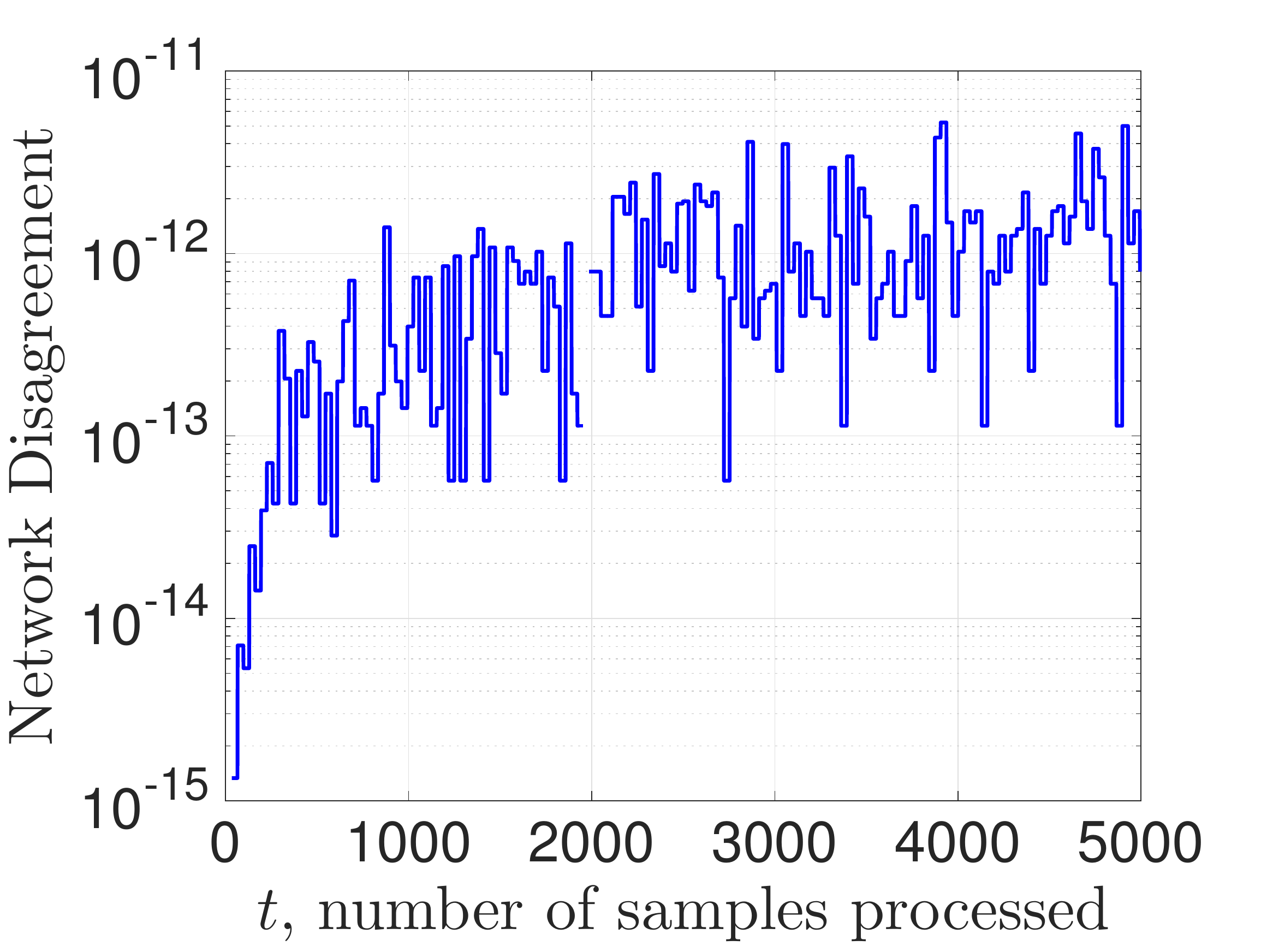}%
\caption{Disagreement vs. samples processed}%
\label{subfigb_sample_path}%
\end{subfigure}\hspace{.5mm}
\begin{subfigure}{0.33\columnwidth}
\includegraphics[width=1.05\linewidth,height = 0.73\linewidth]{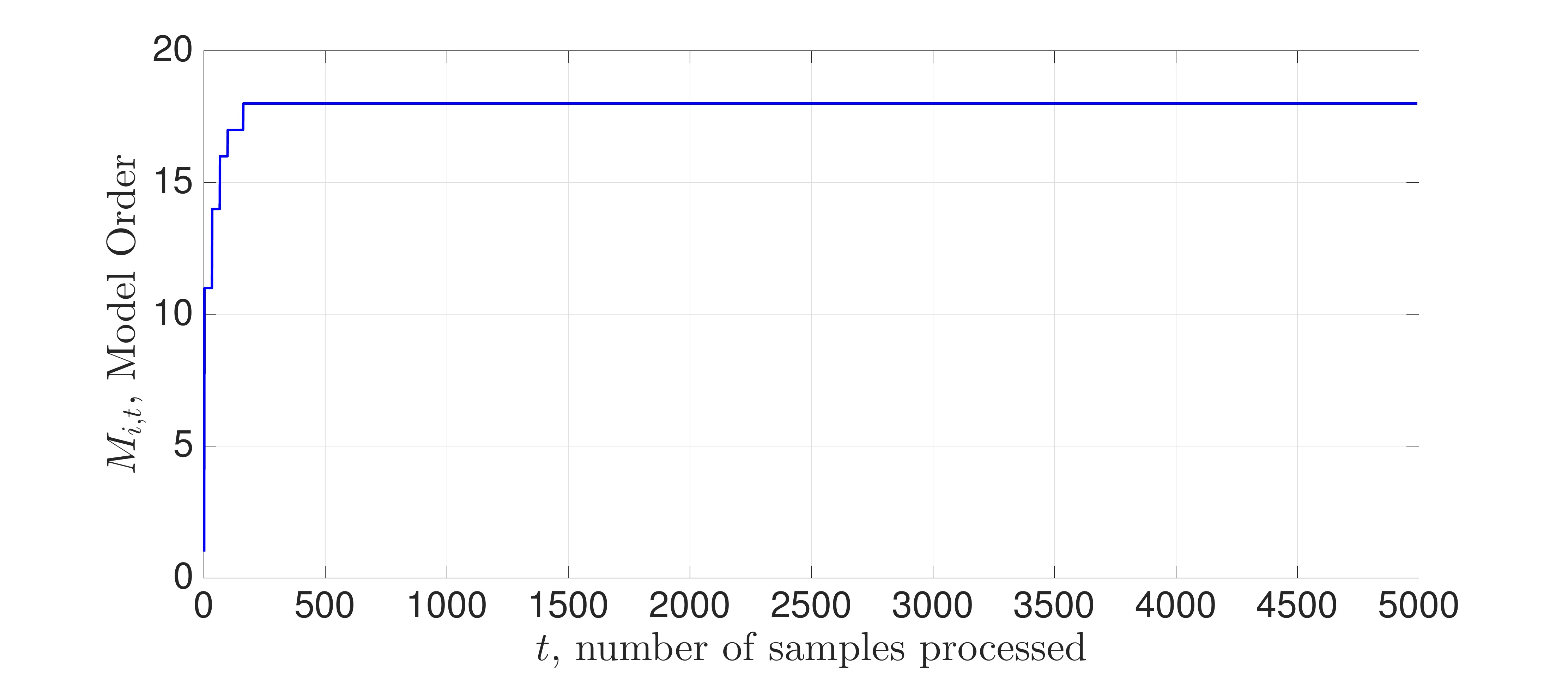}%
\caption{Model Order $M_{i,t}$  vs. samples processed}%
\label{subfigc_sample_path}%
\end{subfigure}%
\caption{ In Fig. \ref{subfiga_sample_path}, we plot the global objective $\sum_{i\in\ccalV}(\mbE_{\bbx_i, y_i}[ \ell_i(f_{i,t}\big(\bbx), y_i\big)])$ versus the number of samples processed, and observe convergence. In Fig. \ref{subfigb_sample_path} we display the Hilbert-norm network disagreement $\sum_{(i,j)\in\ccalE} \|f_{i,t} - f_{j,t}\|_{\ccalH}^2$ with a penalty parameter $c$ that doubles every $200$ samples. As $c$ increases, agents attain consensus. In Fig. \ref{subfigc_sample_path}, we plot the model order of a randomly chosen agent's regression function, which stabilizes to $18$ after $162$ samples. }\label{fig:sample_path}
\end{figure*}

\vspace{-3mm}
\section{Numerical Experiments}\label{sec:simulations}

We consider the task of kernel logistic regression (KLR) (Section \ref{sec_logistic_regression}) from multi-class training data scattered across a multi-agent system in two settings: classification of data from a Gaussian mixture model and texture classification. In Section \ref{sec_svm}, we consider kernel support vector machines (KSVM).\footnote{We thank Garrett Warnell and Ethan Stump of the U.S. Army Research Laboratory for invaluable assistance in the algorithm implementation.} 
%
\vspace{-3mm}
\subsection{Kernel Logistic Regression}\label{sec_logistic_regression}

For KLR, the merit of a particular regressor for agent $i$ is quantified by its contribution to the class-conditional probability. We define a set of class-specific functions $f_{i,k} : \ccalX \rightarrow \reals$, and denote them jointly as $\bbf_i \in \ccalH^D$, where $\{1,\dots,D\}$ denotes the set of classes. Then,
define the probabilistic model
\begin{align}\label{eq:multi_logistic_prob}
P(y_i = d \,|\, \bbx_i) := \frac{\exp(f_{i,d}(\bbx_i))}{\sum_{d^\prime} \exp( f_{i,d^\prime}(\bbx_i))}.
\end{align}
which models the odds ratio of a sample being in class $d$ versus all others. The negative log likelihood defined by \eqref{eq:multi_logistic_prob} is the instantaneous loss (see, e.g., \cite{Murphy2012}) at sample $(\bbx_{i,n},y_{i,n})$:
\begin{align}\label{eq:multi_logistic}
\!\!\!\!\ell_i(\bbf_i,\bbx_{i,n},y_{i,n}) =\! -\! \log P(y_i = y_{i,n} | \bbx_{i,n}). 
\end{align}
%
For a given set of activation functions, classification decisions $\tilde{d}$ for $\bbx_i$ is given by the maximum likelihood estimate, i.e., $\tilde{d}=\argmax_{d\in\{1,\dots,D\}} f_{i,d}(\bbx)$. 

{\bf Gaussian Mixture Model}
Following \cite{Zhu2005,POLK}, we generate a data set from Gaussian mixture models, which consists $N=5000$ feature-label pairs for training and $2500$ for testing.  Each label $y_n$ was drawn uniformly at random from the label set. The corresponding feature vector $\bbx_n \in \reals^p$ was then drawn from a planar ($p=2$), equitably-weighted Gaussian mixture model, i.e., $\bbx \given y \; \sim \; (1/3) \sum_{j=1}^3 \ccalN(\boldsymbol{\mu}_{y,j}, \sigma^2_{y,j}\bbI)$ where $\sigma^2_{y,j}=0.2$ for all values of $y$ and $j$. The means $\boldsymbol{\mu}_{y,j}$ are themselves realizations of their own Gaussian distribution with class-dependent parameters, i.e., $\boldsymbol{\mu}_{y,j} \sim \ccalN( \boldsymbol{\theta}_y, \sigma^2_y\bbI )$, where $\left\{ \bbtheta_1,\ldots,\bbtheta_D\right\}$ are equitably spaced around the unit circle, one for each class label, and $\sigma_y^2=1.0$. We fix the number of classes $D=5$, meaning that the feature distribution has, in total, $15$ distinct modes.  The data is plotted in Figure \ref{ch7_subfiga_gmm}. 

Each agent in a $V=20$ network observes a unique stream of training examples from this common data set. Here the communications graph is a random network with edges generated randomly between nodes with probability $1/5$ repeatedly until we obtain one that is connected, and then symmetrize it. We run Algorithm \ref{alg:soldd} when the entire training set is fed to each agent in a streaming fashion, a Gaussian kernel is used with bandwidth $d=0.6$, with constant learning rate $\eta=3$, compression budget chosen as $\eps=\eta^{3/2}$ with parsimony constant $K=0.04$, mini-batch size $32$, and regularizer $\lambda=10^{-6}$. The penalty coefficient is initialized as $c=0.01$ and doubled after every 200 training examples. 
 
We plot the results of this implementation in Figures \ref{ch7_subfigb_gmm} and \ref{fig:sample_path}. In Figure \ref{subfiga_sample_path}, we plot the global objective $\sum_{i\in\ccalV}(\mbE_{\bbx_i, y_i}[ \ell_i(f_{i,t}\big(\bbx), y_i\big)])$ relative to the number of training examples processed, and observe stable convergence to a global minimum. In Figure \ref{subfigb_sample_path} we display Hilbert-norm network disagreement $\sum_{(i,j)\in\ccalE} \|f_{i,t} - f_{j,t}\|_{\ccalH}^2$ versus observed sample points. Since each regression function is initialized as null, initially the disagreement is trivially null, but it remains small over the function sample path as model training occurs. Moreover, the model order of an arbitrarily chosen agent $i=15$ versus samples processed is given in Figure \ref{subfigc_sample_path}: observe that the model order stabilizes after only a couple hundred training examples to $18$, which is only a couple more than $15$, the number of modes of the joint data density function. The resulting decision surface of node $15$ is given in Figure \ref{ch7_subfigb_gmm}, which achieves $95.2\%$ classification accuracy on the test set which is comparable to existing centralized batch approaches (see Table 2 of \cite{POLK}) to kernel logistic regression.

\begin{figure*}%
\centering \hspace{-4mm}
\begin{subfigure}{0.33\columnwidth}
\includegraphics[width=0.9\linewidth, height = 0.75\linewidth]{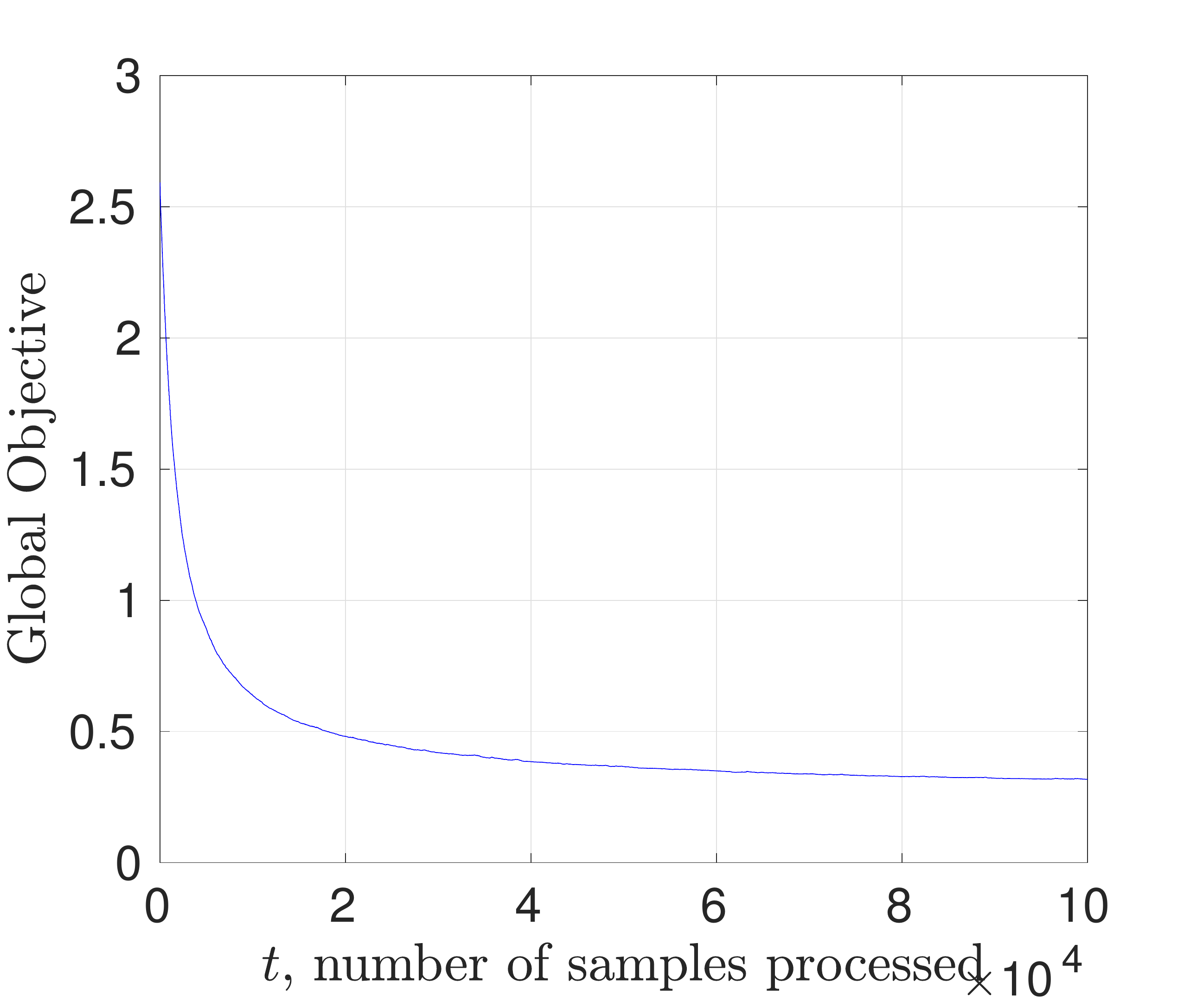}%
\caption{Global Objective vs. samples processed}%
\label{subfiga_loss_brodatz}%
\end{subfigure}\hspace{.5mm}
\begin{subfigure}{0.32\columnwidth}
\includegraphics[width=1.05\linewidth,height = 0.73\linewidth]{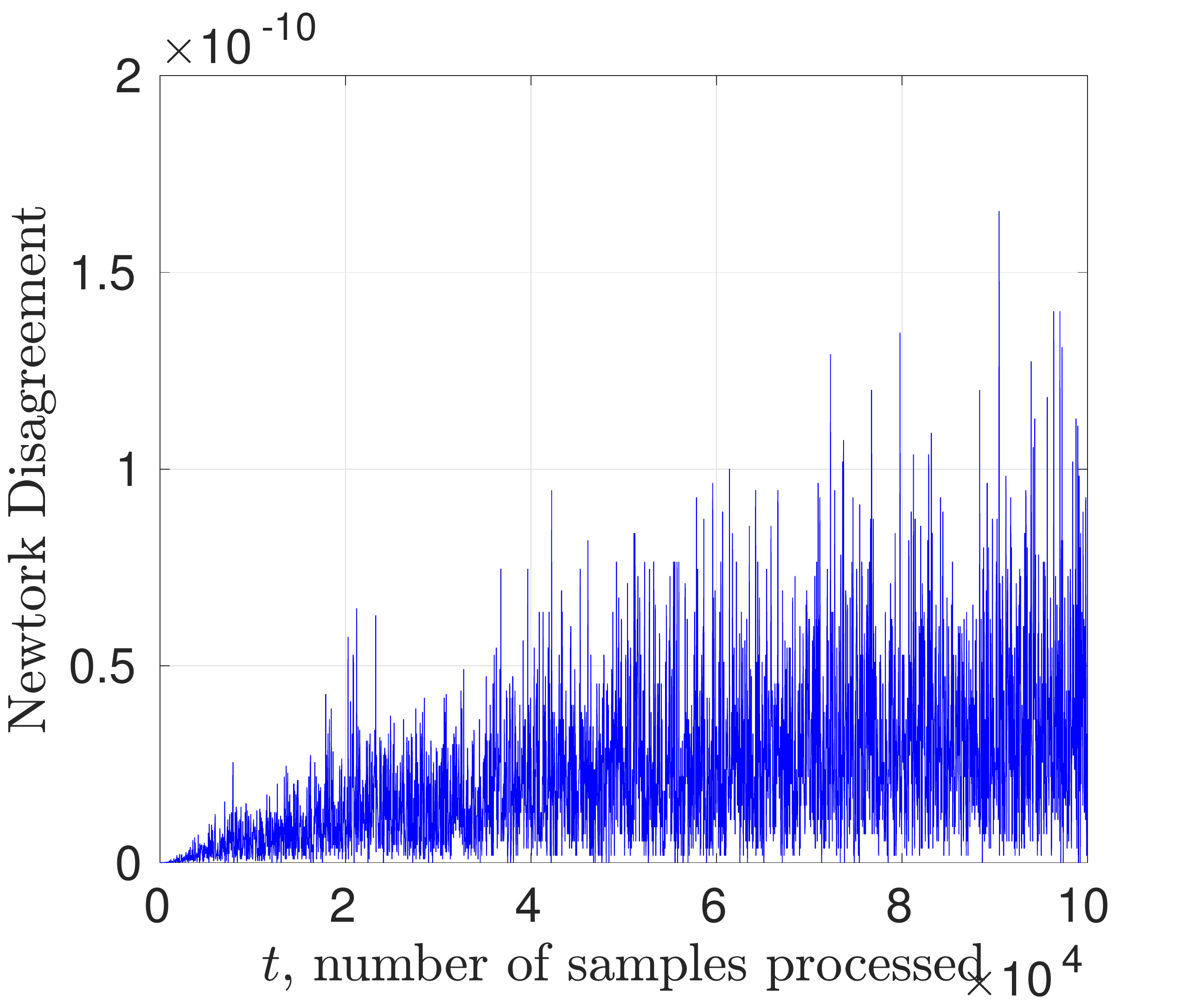}%
\caption{Disagreement vs. samples processed}%
\label{subfigb_disagreement_brodatz}%
\end{subfigure}\hspace{.5mm}
\begin{subfigure}{0.33\columnwidth}
\includegraphics[width=1.05\linewidth,height = 0.75\linewidth]{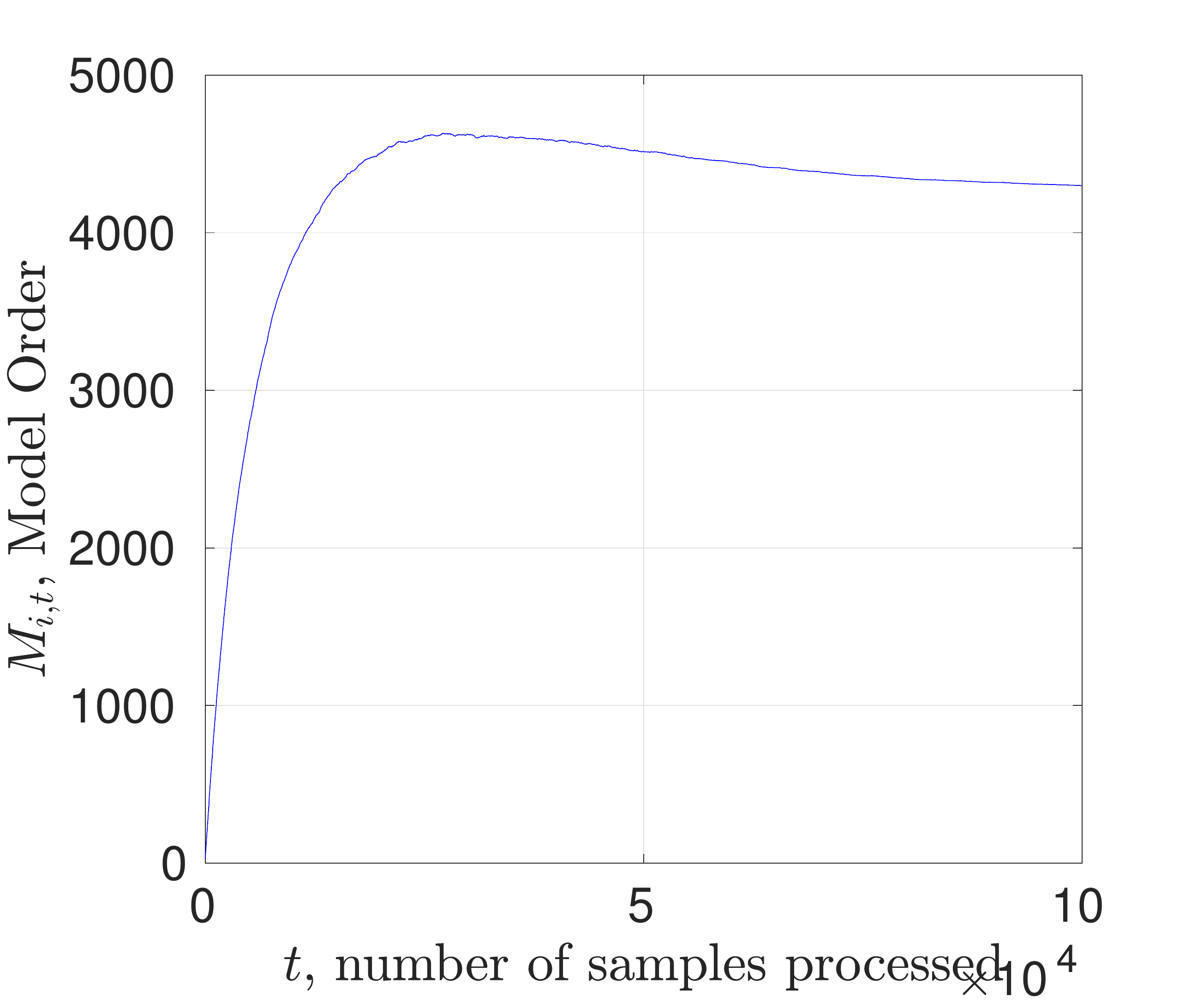}%
\caption{Model Order $M_{i,t}$  vs. samples processed}%
\label{subfigc_model_order_brodatz}%
\end{subfigure}%
\caption{In Fig. \ref{subfiga_loss_brodatz}, we plot the global objective $\sum_{i\in\ccalV}(\mbE_{\bbx_i, y_i}[ \ell_i(f_{i,t}\big(\bbx), y_i\big)])$ versus the number of samples processed, and observe convergence. In Fig. \ref{subfigb_disagreement_brodatz} we display the Hilbert-norm network disagreement $\sum_{(i,j)\in\ccalE} \|f_{i,t} - f_{j,t}\|_{\ccalH}^2$ with a penalty parameter $c=0.02$. In Fig. \ref{subfigc_model_order_brodatz}, we plot the model order of a randomly chosen agent's regression function, which stabilizes to $4299$. }\label{fig_brodatz}
\end{figure*}


{\bf Texture Classification}
We generated the \textit{brodatz} data set using a subset of the images provided in \cite{Brodatz1966}. Specifically, we used 13 texture images (i.e. D=13), and from them generated a set of 256 textons \cite{Leung1999}. Next, for each overlapping patch of size 24-pixels-by-24-pixels within these images, we took the feature to be the associated $p=256$-dimensional texton histogram. The corresponding label was given by the index of the image from which the patch was selected. We then randomly selected $N=10000$ feature-label pairs for training and $5000$ for testing. 
Each agent in network with $V=5$ observes a unique stream of training examples from this common data set. Here the communication graph is a random network with edges generated randomly between nodes with probability $1/5$ repeatedly until we obtain one that is connected, and then symmetrize it. To train the classifier we run Algorithm \ref{alg:soldd} ten epoches: in each epoch we fed the entire training set to each agent in a streaming fashion. A Gaussian kernel is used with bandwith $\sigma^2 = 0.1$, with constant learning rate $\eta =4$, compression budget $\epsilon = \eta^{3/2}$ with parsimony constant $K = 0.04$, mini-batch size $32$ and regularizer $\lambda = 10^{-5}$. The penalty coefficient is set to $c=0.02$.

We plot the results of this experiment in Figure \ref{fig_brodatz}. In Figure \ref{subfiga_loss_brodatz} we display the global objective $\sum_{i\in\ccalV}(\mbE_{\bbx_i, y_i}[ \ell_i(f_{i,t}\big(\bbx), y_i\big)])$ relative to the number of observed examples, and observe convergence to a global minimum. In Figure \ref{subfigb_disagreement_brodatz} we plot the Hilbert norm network disagreement $\sum_{(i,j)\in\ccalE} \|f_{i,t} - f_{j,t}\|_{\ccalH}^2$. Since the initial regression function is null for all agents the disagreement is zero and as observed in Figure \ref{subfigb_disagreement_brodatz} it remains small over the training. Moreover, the model order of an agent chosen at random versus samples processed is given in Figure \ref{subfigc_model_order_brodatz}. The resulting decission function achives $93.5\%$ classification accuracy over the test set which is comparable with the accuracy of the centralized version ($95.6\%$) \cite{POLK}. However the model order requiered is more than twice the model order in the centralized case (4358 in average v.s. 1833\cite{POLK}). Compared to other distributed classification algorithms the current algorithm outperforms them. For instance D4L achieves around $75\%$ classification accuracy \cite{KoppelEtal16a}.

\subsection{Kernel Support Vector Machines}\label{sec_svm}
Now we address the problem of training a multi-class kernel support vector machine online in a multi-agent systems. The merit of a particular regressor is defined by its ability to maximize its classification margin, which may be formulated by first defining a set of class-specific activation functions $f_{i,d} : \ccalX \rightarrow \reals$, and denote them jointly as $\bbf_i \in \ccalH^D$.  In Multi-KSVM, points are assigned the class label of the activation function that yields the maximum response. KSVM is trained by taking the instantaneous loss $\ell$ to be the multi-class hinge function which defines the margin separating hyperplane in the kernelized feature space, i.e.,
\begin{align}\label{ch6_eq:multi_svm}
  \ell_i(\bbf_i,\bbx_n,y_n) &= \max(0,1+f_{i,r}(\bbx_n)-f_{i,y_n}(\bbx_n)),  
\end{align}
where $r = \argmax_{d'\neq y} f_{i,d'}(\bbx)$. See \cite{Murphy2012} for further details.

We consider an implementation where each agent in a $V=20$ network observes a unique stream of training examples from the Gaussian mixtures data set (see Figure \ref{ch7_subfiga_gmm}). Moreover, the communications graph is fixed as a random network with edges generated randomly between nodes with probability $1/5$ repeatedly until we obtain one that is connected, and then symmetrize it. We run Algorithm \ref{alg:soldd} when the entire training set is fed to each agent in a streaming fashion, a Gaussian kernel is used with bandwidth $\tilde{\sigma}^2=0.6$, with constant learning rate $\eta=3$, compression budget chosen as $\eps=\eta^{3/2}$ with parsimony constant $K=0.04$, mini-batch size $32$, and regularizer $\lambda=10^{-6}$. The penalty coefficient is initialized as $c=0.01$ and doubled after every 200 training examples.

\begin{figure*}%
\centering \hspace{-4mm}
\begin{subfigure}{0.33\columnwidth}
\includegraphics[width=1.05\linewidth, height = 0.75\linewidth]{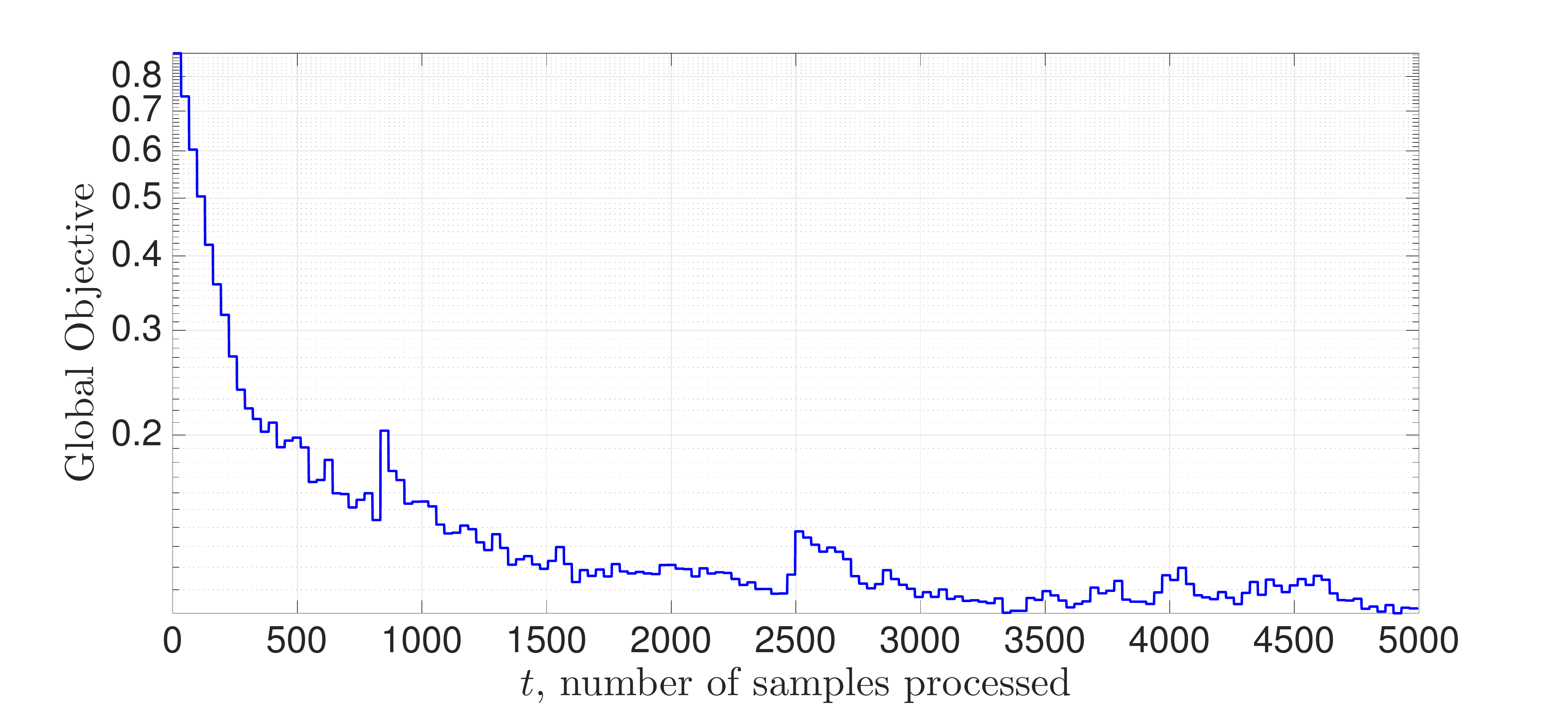}%
\caption{Global objective vs. samples}%
\label{ch7_subfiga_sample_path_svm}%
\end{subfigure}\hspace{.5mm}
\begin{subfigure}{0.32\columnwidth}
\includegraphics[width=1.05\linewidth,height = 0.73\linewidth]{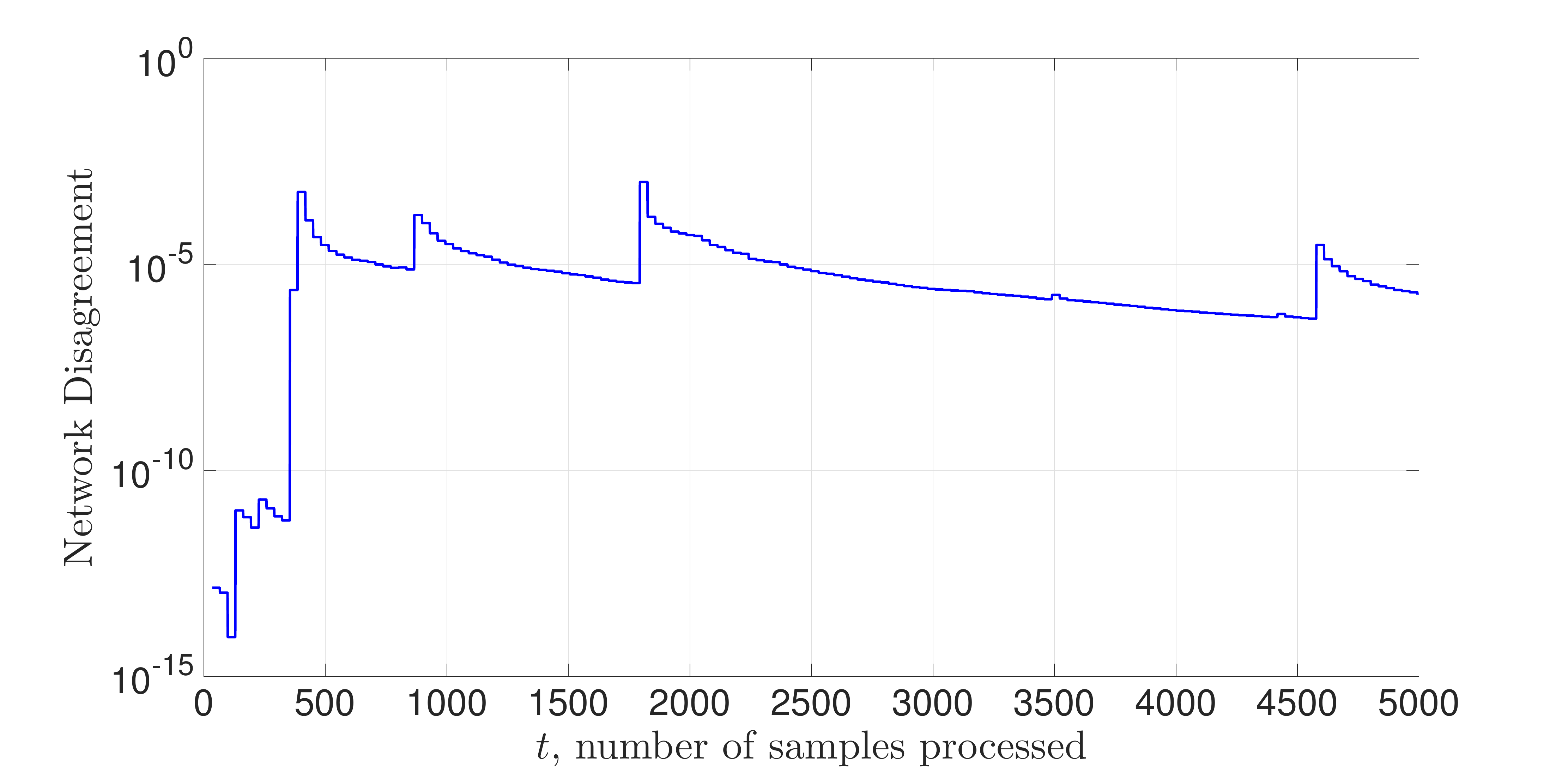}%
\caption{Disagreement vs. samples}%
\label{ch7_subfigb_sample_path_svm}%
\end{subfigure}\hspace{.5mm}
\begin{subfigure}{0.33\columnwidth}
\includegraphics[width=1.05\linewidth,height = 0.75\linewidth]{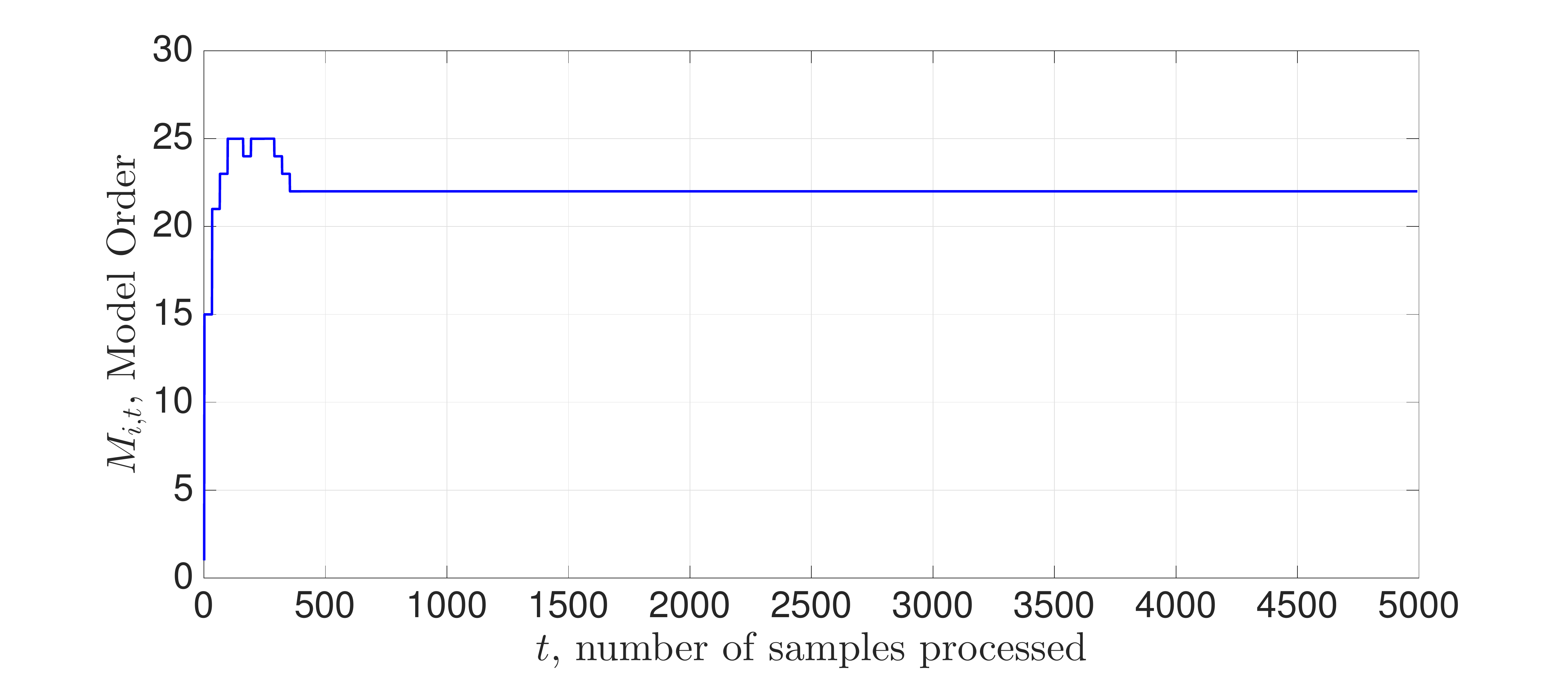}%
\caption{Model Order $M_{i,t}$  vs. samples}%
\label{ch7_subfigc_sample_path_svm}%
\end{subfigure}%
\caption[Stable decentralized memory-efficient kernelized SVM.]{ In Fig. \ref{ch7_subfiga_sample_path_svm}, we plot the global objective $\sum_{i\in\ccalV}(\mbE_{\bbx_i, y_i}[ \ell_i(f_{i,t}\big(\bbx), y_i\big)]) $ versus the number of samples processed, and observe convergence, albeit more noisily than for the differentiable logistic loss. In Fig. \ref{ch7_subfigb_sample_path_svm} we display the Hilbert-norm network disagreement $\sum_{(i,j)\in\ccalE} \|f_{i,t} - f_{j,t}\|_{\ccalH}^2$ with a penalty parameter $c$ that doubles every $200$ samples. As $c$ increases, agents attain consensus with respect to the Hilbert norm. In Fig. \ref{ch7_subfigc_sample_path_svm}, we plot the model order of a randomly chosen agent's regression function, which stabilizes to $22$ after $354$ samples. Here we obtain a slightly higher complexity classifier that achieves slightly better accuracy.}\label{ch7_fig:sample_path_svm}
\end{figure*}

We plot the results of this implementation in Figures \ref{ch7_subfigc_gmm} and \ref{ch7_fig:sample_path_svm}. In Figure \ref{ch7_subfiga_sample_path_svm}, we observe that the global objective $\sum_{i\in\ccalV}(\mbE_{\bbx_i, y_i}[ \ell_i(f_{i,t}\big(\bbx), y_i\big)])$ converges stably to a global minimum as the number of samples processed increases. In Figure \ref{ch7_subfigb_sample_path_svm} we display Hilbert-norm network disagreement $\sum_{(i,j)\in\ccalE} \|f_{i,t} - f_{j,t}\|_{\ccalH}^2$ versus observed sample points. Since each regression function is initialized as null, initially the disagreement is trivially null, but it remains small over the function sample path as model training occurs, and periodically spikes when the penalty parameter is increased. Moreover, the model order of an arbitrarily chosen agent $i=6$ versus samples processed is given in Figure \ref{ch7_subfigc_sample_path_svm}: the model order stabilizes after only a couple hundred training examples to $22$, which is only a couple more than $15$, the number of modes of the joint data density function. The resulting decision surface of node $6$ is given in Figure \ref{ch7_subfigc_gmm}, which achieves $95.7\%$ classification accuracy, which is approximately state of the art.

\section{Conclusion} \label{sec:conclusion}
In this paper, we extended the ideas in \cite{POLK} to multi-agent settings with the intent of developing a method such that a network of autonomous agents, based on their local data stream, may learn a kernelized statistical model which is optimal with respect to information aggregated across the entire network. To do so, we proposed an unusual penalty function whose structure is amenable to efficient parameterizations when developing stochastic approximation-based updates. By applying functional stochastic gradient method to this node-separable penalty combined with greedily constructed subspace projections, we obtain a decentralized online algorithm for memory-efficient nonparametric function approximation that is globally convergent. We obtain a controllable trade-off between optimality and memory requirements through the design of the greedy subspace projections. Moreover, for large penalty parameter selections, agents achieve consensus.

The empirical performance of this protocol, the Greedy Projected Penalty Method, yields state of the art statistical accuracy for a team of interconnected agents learning from streaming data for both multi-class kernel logistic regression and multi-class kernel support vector machines problems. These results provide a mathematical and empirical foundation for accurate and stable multi-agent statistical inference in online settings while preserving memory-efficiency.  


\vspace{-2mm}
\section*{Appendix A: Details of Matching Pursuit} \label{apx_matching_pursuit}

The removal procedure is as follows: 
at each step, a single dictionary element $j$ of $\bbD$ is selected to be removed which contributes the least to the Hilbert-norm error $\min_{f\in\ccalH_{\bbD\setminus \{j\}}}\|\tilde{f} - f \|_{\ccalH}$ of the original function $\tilde{f}$, when dictionary $\bbD$ is used. Since at each stage the kernel dictionary is fixed, this amounts to a computation involving weights $\bbw \in \reals^{M-1}$ only; that is, the error of removing dictionary point $\bbd_j$ is computed for each $j$ as 
 $\gamma_j =\min_{\bbw_{\ccalI \setminus \{j\}}\in\reals^{{M}-1}} \|\tilde{f}(\cdot) - \sum_{k \in \ccalI \setminus \{j\}} w_k \kappa(\bbd_k, \cdot) \|.$
  We use the notation $\bbw_{\ccalI \setminus \{j\}}$ to denote the entries of $\bbw\in \reals^M$ restricted to the sub-vector associated with indices $\ccalI \setminus \{j\}$. Then, we define the dictionary element which contributes the least to the approximation error as $j^*=\argmin_j \gamma_j$. If the error incurred by removing this kernel dictionary element exceeds the given compression budget $\gamma_{j^*}>\eps_t$, the algorithm terminates. Otherwise, this dictionary element $\bbd_{j^*}$ is removed, the weights $\bbw$ are revised based on the pruned dictionary as $\bbw = \argmin_{\bbw \in \reals^{{M}}} \lVert \tilde{f}(\cdot) - \bbw^T\boldsymbol{\kappa}_{\bbD}(\cdot) \rVert_{\ccalH}$, and the process repeats as long as the current function approximation is defined by a nonempty dictionary. This procedure is summarized in Algorithm \ref{alg:komp}.

\begin{algorithm}[t]
\caption{Kernel Orthogonal Matching Pursuit (KOMP) \hspace{-2mm}}
\begin{algorithmic}
\label{alg:komp}
\REQUIRE  function $\tilde{f}$ defined by dict. $\tbD \in \reals^{p \times \tilde{M}}$, coeffs. $\tbw \in \reals^{\tilde{M}}$, approx. budget  $\epsilon_t > 0$ \\
\STATE \textbf{initialize} $f=\tilde{f}$, dictionary $\bbD = \tbD$ with indices $\ccalI$, model order $M=\tilde{M}$, coeffs.  $\bbw = \tbw$.
\WHILE{candidate dictionary is non-empty $\ccalI \neq \emptyset$}
{\FOR {$j=1,\dots,\tilde{M}$}
	\STATE Find minimal approximation error with dictionary element $\bbd_j$ removed \vspace{-2mm}
	$$\gamma_j = \min_{\bbw_{\ccalI \setminus \{j\}}\in\reals^{{M}-1}} \|\tilde{f}(\cdot) - \sum_{k \in \ccalI \setminus \{j\}} w_k \kappa(\bbd_k, \cdot) \|_{\ccalH} \; .$$ \vspace{-5mm}
\ENDFOR}
	\STATE Find index minimizing approx. error: $j^* = \argmin_{j \in \ccalI} \gamma_j$ 
	\INDSTATE{{\bf{if }} minimal approx. error exceeds threshold $\gamma_{j^*}> \epsilon_t$}
	\INDSTATED{\bf stop} 
	\INDSTATE{\bf else} 
	
	\INDSTATED Prune dictionary $\bbD\leftarrow\bbD_{\ccalI \setminus \{j^*\}}$
	\INDSTATED Revise set $\ccalI \leftarrow \ccalI \setminus \{j^*\}$, model order ${M} \leftarrow {M}-1$.
	\INDSTATED Update weights $\bbw$ defined by current dictionary $\bbD$
	$$\bbw = \argmin_{\bbw \in \reals^{{M}}} \lVert \tilde{f}(\cdot) - \bbw^T\boldsymbol{\kappa}_{\bbD}(\cdot) \rVert_{\ccalH}$$\vspace{-5mm}
	\INDSTATE {\bf end}
\ENDWHILE	
\RETURN ${f},\bbD,\bbw$ of model order $M \leq \tilde{M}$ such that $\|f - \tilde{f} \|_{\ccalH}\leq \eps_t$
\end{algorithmic}
\end{algorithm}

\vspace{-2mm}
\section*{Appendix B: Proof of Proposition \ref{prop_projection}}\label{apx_prop_projection}
%
Consider the square-Hilbert-norm difference of the stacked projected stochastic gradient $ \tilde{\nabla}_{f} \hat{\psi}_{c} (f_{t}(\bbx_{t}\!),y_{t})$  and its un-projected variant $\nabla_f \hat{\psi}_c(f_t(\bbx_t), \bby_t) $ defined in \eqref{eq:proj_grad} and \eqref{eq:penalty_function_stoch_grad}, respectively,
\begin{align}\label{eq:norm_stoch_grad_diff}
 \| &\tilde{\nabla}_{f} \hat{\psi}_{c} (f_{t}(\bbx_{t}),y_{t})-\nabla_f \hat{\psi}_c(f(\bbx_t), \bby_t)  \|_{\ccalH}^2  \\
 &= \Big\| \text{vec}\!\Big(\! f_{i,t} \!-\! \ccalP_{ \ccalH_{\bbD_{i,t+1}}}\! \!\Big[\!
 f_{i,t} 
\!-\! {\eta}_t \!{\nabla}_{f_i}\hat{\psi}_{i,c}(f_{i,t}(\bbx_{i,t}\!),\!y_{i,t}\!) \!\Big]\!\Big)\!/\eta_t \nonumber \\
&\qquad - \text{vec}\left(\nabla_{f_i} \hat{\psi}_{i,c}(f_{i,t}(\bbx_{i,t}), y_{i,t})\right)  \Big\|_{\ccalH}^2 \nonumber \\
&  \leq \! V^2\! \max_{i\in\ccalV}  \Big\| \!\Big(\! f_{i,t} \!-\! \ccalP_{ \ccalH_{\bbD_{i,t+1}}}\! \!\Big[\!
 f_{i,t} 
\!-\! {\eta}_t \!{\nabla}_{f_i}\hat{\psi}_{i,c}(f_{i,t}(\bbx_{i,t}\!),\!y_{i,t}\!) \!\Big]\!\Big)\!/\eta_t \nonumber \\
&\qquad - \nabla_{f_i} \hat{\psi}_{i,c}(f_{i,t}(\bbx_{i,t}), y_{i,t})  \Big\|_{\ccalH}^2 \nonumber
\end{align}
where we apply the fact that the functional gradient is a concatenation of functional gradients associated with each agent in \eqref{eq:norm_stoch_grad_diff} for the first equality, and for the second inequality we consider the worst-case estimate across the network. Now, let's focus on the term inside the Hilbert-norm on the right-hand side. Multiply and divide $\nabla_{f_i} \hat{\psi}_{i,c}(f_{i,t}(\bbx_{i,t}), y_{i,t}) $, the last term, by $\eta_t$, and reorder terms to write
\begin{align}\label{eq:norm_stoch_grad_expand}
 \Big\| \!\Big(\! f_{i,t} - &\ccalP_{ \ccalH_{\bbD_{i,t+1}}}\! \!\Big[\!
 f_{i,t} 
\!-\! {\eta}_t \!{\nabla}_{f_i}\hat{\psi}_{i,c}(f_{i,t}(\bbx_{i,t}\!),\!y_{i,t}\!) \!\Big]\!\Big)\!/\eta_t \nonumber \\
& - \nabla_{f_i} \hat{\psi}_{i,c}(f_{i,t}(\bbx_{i,t}), y_{i,t})  \Big\|_{\ccalH}^2 \nonumber \\
& = \Big\|\! \frac{1}{\eta_t}\!\!\left(\! f_{i,t}\! -\!\eta_t \nabla_{f_i} \hat{\psi}_{i,c}(f_{i,t}(\bbx_{i,t}), y_{i,t})\!  \right)\!
\nonumber \\
&\qquad -\!\frac{1}{\eta_t}\ccalP_{ \ccalH_{\bbD_{i,t+1}}}\! \Big[  f_{i,t} \!-\! \eta_t \hat{\psi}_{i,c}(f_{i,t}(\bbx_{i,t}), y_{i,t})  \Big]  \Big\|_{\ccalH}^2 \nonumber \\
&\qquad \qquad=\frac{1}{\eta_t^2}\| \tilde{f}_{i,t+1} - f_{i,t+1} \|_{\ccalH}^2 
  \end{align}
  where we have substituted the definition of $\tilde{f}_{i,t+1}$ and $f_{i,t+1}$ in  \eqref{eq:sgd_tilde} and \eqref{eq:projection_hat}, respectively, and pulled the nonnegative scalar $\eta_t$ outside the norm. Now, observe that the KOMP residual stopping criterion in Algorithm \ref{alg:komp} is $\lVert \tilde{f}_{i,t+1} - f_{i,t+1} \rVert_{\ccalH} \leq \epsilon_t$, which we may apply to the last term on the right-hand side of \eqref{eq:norm_stoch_grad_expand}. This result with the inequality \eqref{eq:norm_stoch_grad_diff} yields \eqref{eq:prop_projection}.$\hfill \blacksquare$\vspace{-2mm}

%

\section*{Appendix C: Proof of Lemma \ref{lemma1}}\label{apx_lemma_descent}
Begin by considering the square of the Hilbert-norm difference between $f_{t+1}$ and $f^*_c=\argmin \psi_c(f)$ which minimizes \eqref{eq:penalty_function}, and expand the square to write
\begin{align}\label{eq:iterate_square_expand}
\| f_{t+1} - f^*_c \|_{\ccalH}^2 
	 &=\| f_t - \eta_t \tilde{\nabla}_{f} \hat{\psi}_{c} (f_{t}(\bbx_{t}),\bby_{t}) \|_{\ccalH}^2 \nonumber \\
	& =\!\|  f_t \! -\! f^* \|_{\ccalH}^2 \!-\! 2 \eta_t \langle f_t \!-\! f^*_c, \tilde{\nabla}_{f} \hat{\psi}_{c} (f_{t}(\!\bbx_{t}),\bby_{t}\!) \rangle_{\ccalH} \nonumber \\
	& \quad+ \eta_t^2 \| \tilde{\nabla}_{f} \hat{\psi}_{c} (f_{t}(\!\bbx_{t}),\bby_{t}) \|_{\ccalH}^2
\end{align}
Add and subtract the functional stochastic gradient of the penalty function  ${\nabla}_{f} \hat{\psi}_{c} (f_{t}(\!\bbx_{t}),\bby_{t}\!)$ defined in \eqref{eq:penalty_function_stoch_grad} to the second term on the right-hand side of \eqref{eq:iterate_square_expand} to obtain
\begin{align}\label{eq:iterate_stoch_grad}
\| f_{t+1} - f^*_c \|_{\ccalH}^2 
	 &=\!\| f_t \!-\! f^*_c \|_{\ccalH}^2
	 \!- \! 2 \eta_t \langle f_t\! -\! f^*_c,\! \!{\nabla}_{f} \hat{\psi}_{c} (f_{t}(\!\bbx_{t}),\bby_{t}\!)  \rangle_{\ccalH} \nonumber\\
	 &\!\! \!\!\!\!\!\!-\! 2 \eta_t \!\langle f_t \!-\!\! f^*_c, \!\!\tilde{\nabla}_{f} \hat{\psi}_{c} (f_{t}(\!\bbx_{t}),\bby_{t}\!) \!-\!\! {\nabla}_f \hat{\psi}_c(f_{t}(\!\bbx_t),\bby_t\!)\! \rangle_{\ccalH} \nonumber\\
	&\quad+ \eta_t^2 \| \tilde{\nabla}_f\ell(f_{t}(\bbx_t),\bby_t) \|_{\ccalH}^2 
\end{align}
We deal with the third term on the right-hand side of \eqref{eq:iterate_stoch_grad}, which represents the directional error associated with the sparse stochastic projections, by applying the Cauchy-Schwartz inequality together with Proposition \ref{prop_projection} to obtain
\begin{align}\label{eq:iterate_prop1_subst}
\| f_{t+1} \!-\! f^*_c \|_{\ccalH}^2 
	 &\!=\! \| f_t \!-\! f^*_c \|_{\ccalH}^2
	 \!-\! 2 \eta_t \langle f_t \!-\! f^*_c,\!\! {\nabla}_{f} \hat{\psi}_{c} (f_{t}(\!\bbx_{t}),\bby_{t}\!)  \rangle_{\ccalH} \nonumber\\
	 &\!\!\!\!+\!  2 \eps_t V \| f_t \!-\! \!f^*_c\|_{\ccalH}\! +\! \eta_t^2 \| \tilde{\nabla}_f\ell(f_{t}(\bbx_t),y_t)\! \|_{\ccalH}^2 
\end{align}
Now compute the expectation of \eqref{eq:iterate_prop1_subst} conditional on the algorithm history $\ccalF_t$
\begin{align}\label{eq:expectation_square_expand}
\E{\! \| f_{t+1} \!- \! f^*_c \|_{\ccalH}^2 \given \! \ccalF_t}
 &=\| \! f_t \! - \! f^*_c \|_{\ccalH}^2+  2 \eps_t V\| f_t - f^*_c\|_{\ccalH} + \eta_t^2\sigma^2 \nonumber\\
	&\  \ - \! 2 \eta_t \! \langle f_t \!- \!f^*_c\! , \nabla_f {\psi}_c(f_t) \rangle_{\ccalH} 
\end{align}
where we have applied the fact that the stochastic functional gradient in \eqref{eq:penalty_function_stoch_grad} is an unbiased estimator [cf. \eqref{eq:unbiased}] for the functional gradient of the penalty function in \eqref{eq:penalty_function}, as well as the fact that the variance of the functional projected stochastic gradient is finite stated in \eqref{eq:stochastic_grad_var} (Assumption \ref{as:last}). 
Observe that since $\psi_c(f)$ is an expectation of a convex function, it is also convex, which allows us to write
\begin{equation}\label{eq:convexity}
\psi_c(f_t) -\psi_c(f^*_c) \leq \langle f_t - f^*_c, {\nabla}_f \psi_c(f_{t}) \rangle_{\ccalH}  \;,
\end{equation}
which we substitute into the second term on the right-hand side of the relation given in \eqref{eq:expectation_square_expand} to obtain
\begin{align}\label{eq:expectation_convexity}
\E{\| f_{t+1} - f^*_c \|_{\ccalH}^2 \given \ccalF_t}
	 &\leq\| f_t \!-\! f^*_c \|_{\ccalH}^2
	 \!-\! 2 \eta_t [\psi_c(f_t) \!-\!\psi_c(f^*_c) ]  \nonumber \\
	&\quad + 2 \eps_t V \| f_t - f^*_c \|_{\ccalH} + \eta_t^2 \sigma^2 \; . 
\end{align}
Thus the claim in Lemma \ref{lemma1} is valid.$\hfill \blacksquare$
%

\section*{Appendix D: Proof of Theorem \ref{corollary:convergence_constant}}\label{label_constant_corollary}

The use of the regularizer $(\lambda/2)\|f \|^2_\ccalH$ in \eqref{eq:penalty_function} implies that the penalty is $\lambda$-strongly convex in $f\in\ccalH$, yielding
\begin{equation}\label{eq:strong_cvx}
 \frac{\lambda}{2} \|f_t - f^*_c \|_{\ccalH}^2 \leq \psi_c(f_t) - \psi_c(f^*_c)
\end{equation}
Substituting the relation \eqref{eq:strong_cvx} into the second term on the right-hand side of the expected descent relation stated in Lemma \ref{lemma1}, with constant step-size $\eta_t=\eta$ and budget $\eps_t=\eps$, yields
\begin{align}\label{eq:expected_descent}
\mathbb{E}[\| & f_{t+1} -  f^*_c \|_{\ccalH}^2 \given \ccalF_t]  \\
	 &\qquad \leq (1 - \eta \lambda )\| f_t - f^*_c \|_{\ccalH}^2
	 + 2 \eps V  \| f_t - f^*_c \|_{\ccalH} 
	+ \eta^2 \sigma^2 \; . \nonumber
\end{align}
The expression in \eqref{eq:expected_descent} may be used to construct a stopping stochastic process , which tracks the suboptimality of $\| f_t - f_c^* \|_{\ccalH}^2$ until it reaches a specific threshold, as in the proof of Theorem 2 of \cite{POLK}. In doing so, we obtain convergence to a neighborhood. We may define a stochastic process $\delta_t$ that qualifies as a supermartingale, i.e. $\E{ \delta_{t+1} \given \ccalF_t } \leq \delta_t$ by considering \eqref{eq:expected_descent} and solving for the appropriate threshold by analyzing when the following holds true
\begin{align}\label{eq:martingale_construct1}
\mathbb{E}[&\| f_{t+1} - f^*_c \|_{\ccalH}^2 \given \ccalF_t] \\
	 &\leq (1 - \eta \lambda )\| f_t - f^*_c \|_{\ccalH}^2
	 + 2 \eps V \| f_t - f^*_c \|_{\ccalH} 
	+ \eta^2 \sigma^2 \nonumber \\
	& \leq \| f_t - f^*_c \|^2_{\ccalH} \; . \nonumber
\end{align}
which may be rearranged to obtain the sufficient condition
\begin{align}\label{eq:martingale_construct2}
  - \eta \lambda \| f_t - f^*_c \|_{\ccalH}^2
	 + 2 \eps V \| f_t - f^*_c \|_{\ccalH} 
	+ \eta^2 \sigma^2 \leq 0\; . 
\end{align}
Note that \eqref{eq:martingale_construct2} defines a quadratic polynomial in $\| f_t - f^*_c \|_{\ccalH}$, which, using the quadratic formula, has roots
\begin{align}\label{eq:quadratic_formula}
\| f_t - f^*_c \|_{\ccalH}=  \frac{\eps V \pm \sqrt{\eps^2 V^2 + \lambda \eta^3  \sigma^2 }}{\lambda \eta} 
\end{align}
Observe \eqref{eq:martingale_construct2} is a downward-opening polynomial in $\|f_t - f^*_c \|_{\ccalH}$ which is nonnegative. Thus, focus on the positive root, substituting the approximation budget selection $\eps=K \eta^{3/2}$ to define the radius of convergence as 
\begin{align}\label{eq:radius}
\!\!\!\Delta:= \!\frac{\eps V\! +\! \sqrt{\eps^2V^2 \!+\! \lambda \eta^3  \sigma^2 }}{\lambda \eta} \!= \!\frac{\sqrt{\eta}}{\lambda}\!\Big(\!\!K V \!\!+\!  \sqrt{\!K^2 V^2 \!+\! \lambda  \sigma^2 }\Big)
\end{align}
 \eqref{eq:radius} allows us to construct a stopping process: define $\delta_t$ as
\begin{align}\label{eq:delta}
\delta_t &= \|f_t - f^*_c \|_\ccalH \\
 &\! \! \! \! \! \times \mathbbm{1} \Big\{\min_{u\leq t}   - \eta \lambda \| f_u - f^*_c \|_{\ccalH}^2
	 + 2 \eps V \| f_u - f^*_c \|_{\ccalH} 
	+ \eta^2 \sigma^2    > \Delta \Big\} \nonumber
\end{align}
where $\mathbbm{1}\{E \}$ denotes the indicator process of event $E \in \ccalF_t$. Note that $\delta_t\geq 0$ for all $t$, since both $\| f_t - f^* \|_{\ccalH}$ and the indicator function are nonnegative. The rest of the proof applies the same reasoning as that of Theorem 2 in \cite{POLK}: in particular, given the definition \eqref{eq:delta}, either 
$\min_{u\leq t}   - \eta \lambda \| f_u - f_c^* \|_{\ccalH}^2
	 + 2 \eps V \| f_u - f^*_c \|_{\ccalH} 
	+ \eta^2 \sigma^2    > \Delta$ holds, in which case we may compute the square root of the condition  in \eqref{eq:martingale_construct1} to write
\begin{align}\label{eq:delta_martingale}
\mathbb{E}[\delta_{t+1} \given \ccalF_t ] \leq \delta_t
\end{align}
Alternatively, $\min_{u\leq t}   - \eta \lambda \| f_u - f^*_c \|_{\ccalH}^2
	 + 2 \eps V \| f_u - f^*_c \|_{\ccalH} 
	+ \eta^2 \sigma^2    \leq \Delta$, in which case the indicator function is null for all $s\geq t$ from the use of the minimum inside the indicator in \eqref{eq:delta}. Thus in either case, \eqref{eq:delta_martingale} is valid, implying $\delta_t$ converges almost surely to null, which, as a consequence we obtain
	%
the fact that either %
$\lim_{t\rightarrow \infty }\|f_t - f^*_c \|_\ccalH - \Delta = 0$ or
the indicator function is null for large $t$, i.e. 
$\lim_{t\rightarrow \infty} \mathbbm{1} \{\min_{u\leq t}   - \eta \lambda \| f_u - f^*_c \|_{\ccalH}^2
	 + 2 \eps V \| f_u - f^*_c \|_{\ccalH} 
	+ \eta^2 \sigma^2    > \Delta \} = 0$ almost surely.
Therefore, we obtain that 
\begin{align}\label{eq:liminf_constant}
\!\!\liminf_{t\rightarrow \infty} \|f_t \!- \!f^*_c \|_\ccalH \!\leq \! \Delta\! = \!\frac{\sqrt{\eta}}{\lambda}\!\Big(K V \!\!+\!  \sqrt{K^2 \!+\! \lambda  \sigma^2 }\Big)\! \ \text{ a.s. } ,
\end{align}
 as stated in Theorem \ref{corollary:convergence_constant}.$\hfill \blacksquare$\vspace{-2mm}


\section*{Appendix E: Proof of Proposition \ref{prop_constraint_violation}}\label{appendix_subopt}
Let $f_c^*$ be the minimizer of $\psi_c(f)$ defined in \eqref{eq:penalty_function} and $f^*$ be the solution of the problem \eqref{eq:main_prob}.  Since the former is the minimizer of $\psi_c(f)$ it holds that 
\begin{align}\label{eq:f_c_suboptimality}
 \psi_c(f_c^*) \leq\psi_c(f^*) \! &=  \!\!\sum_{i\in\ccalV}\!\!\Big(\!\mbE_{\bbx_i, \bby_i}\! \Big[ \ell_i(f^*_i\big(\bbx_i), y_i\big)\!\Big]\! \!+\!\frac{\lambda}{2}\|f^*_i \|^2_{\ccalH} \!  \nonumber\\
&\quad\! \! + \! \frac{c}{2}  \!\!\sum_{j\in n_i}\!  \mathbb{E}_{\bbx_i}\left\{ \! [f^*_i(\bbx_i) \! - \! f^*_j(\bbx_i)]^2 \! \right\}\!\Big) \; . 
\end{align}
Where the equality follows from the definition of $\psi_c(f)$ in \eqref{eq:penalty_function}. Since $f^*$ is solution to the problem \eqref{eq:main_prob} it satisfies that $f_i =f_j$ for all $(i,j)\in\ccalE$, thus 
\begin{equation}
\mathbb{E}_{\bbx_i}\left\{ [f^*_i(\bbx_i) -  f^*_j(\bbx_i)]^2 \right\}=0 \; ,
\end{equation}
for all $(i,j)\in\ccalE$.
As a consequence, replacing $\psi_c(f_c^*)$ by its expression in the first equality in \eqref{eq:f_c_suboptimality} and rearranging terms yields a bound the constraint violation of $f^*_c$ as 
\begin{align}\label{eqn_constraint_violation}
\! \frac{1}{2} \sum_{i\in \ccalV} \!\!\sum_{j\in n_i}\!  \mathbb{E}_{\bbx_i}\left\{ [f^*_{c,i}(\bbx_i) \! - \! f^*_{c,j}(\bbx_i)]^2 \right\} \leq\frac{1}{c}\left(R(f^*)-R(f^*_c)\right)\; ,
\end{align}
where $R(f)$ is the global regularized objective in \eqref{eq:kernel_stoch_opt_global}, i.e.,
\begin{equation}
R(f) =\!\!\sum_{i\in\ccalV}\!\!\Big(\!\mbE_{\bbx_i, \bby_i}\! \Big[ \ell_i(f_i\big(\bbx_i), y_i\big)\!\Big]\! \!+\!\frac{\lambda}{2}\|f_i \|^2_{\ccalH}\Big).
\end{equation}
The fact that by definition $p^*=R(f^*)$ yields \eqref{eqn_constraint_bound}.

%
%
%
%
%

\bibliographystyle{IEEEtran}
\bibliography{bibliography}
\newpage

   \end{document}